\documentclass[12pt, a4paper]{amsart}
\usepackage{amscd}
\usepackage{amsmath}
\usepackage{latexsym}
\usepackage{graphicx}
\usepackage{amsfonts}
\usepackage{amssymb}
\nonstopmode
\newcounter{fig}

\theoremstyle{plain}
\newtheorem{Theorem}{Theorem}[section]
\newtheorem{theorem}[Theorem]{Theorem}

\newtheorem{Corollary}[Theorem]{Corollary}
\newtheorem{Proposition}[Theorem]{Proposition}
\newtheorem{Lemma}[Theorem]{Lemma}
\newtheorem{Void}[Theorem]{}
\theoremstyle{definition}

\newtheorem{Remarks}[Theorem]{Remarks}

\theoremstyle{remark}

\def\Changed/{\ifvmode\else\vadjust{\vbox to 0pt{\vskip -\baselineskip\hbox to 0pt{\hss\vrule height 0pt depth 1.2\baselineskip\hskip 1em}\vss}}\fi}
\def\Math#1{\def\MathString{#1}\futurelet\MathDelim\MathChoose}
\def\MathChoose{\ifmmode\let\MathDo\MathString              \else\let\MathDo\MathSkip\fi              \MathDo}
\def\MathSkip{\ifx\MathDelim/\def\MathDo{$\MathString$\EatOne}              \else\def\MathDo{$\MathString$}\fi              \MathDo}
\def\Text#1{\def\TextString{#1}\futurelet\TextDelim\TextSkip}
\def\TextSkip{\ifx\TextDelim/\def\TextDo{\TextString\EatOne}              \else\let\TextDo\TextString\fi              \TextDo}
\def\EatOne#1{}
\def\SkipToEndScan#1\EndScan{}
\def\Scan#1#2#3{\ifx#1#2#3\expandafter\SkipToEndScan\fi\Scan#1}
\def\Upper#1{\Scan#1aAbBcCdDeEfFgGhHiIjJkKlLmMnNoOpPqQrRsStTuUvVwWxXyYzZ#1#1\EndScan}
\def\Phrase#1 #2/#3/#4=#5 #6/#7/#8.{\expandafter\edef\csname#2#3\endcsname{\noexpand\Text{#6#7}}
\expandafter\edef\csname\Upper#2#3\endcsname{\noexpand\Text{\Upper#6#7}}
\expandafter\edef\csname#1#2#3\endcsname{\noexpand\Text{#5 #6#7}}
\expandafter\edef\csname\Upper#1#2#3\endcsname{\noexpand\Text{\Upper#5 #6#7}}
\expandafter\edef\csname#2#4\endcsname{\noexpand\Text{#6#8}}
\expandafter\edef\csname\Upper#2#4\endcsname{\noexpand\Text{\Upper#6#8}}
}

\begin{document}
\title[Hermitian $K$-theory of the integers]{Hermitian $K$-theory of the integers}
\author{A. J. Berrick}
\author{M. Karoubi}
\address{Department of Mathematics, National University of Singapore, Kent Ridge
117543, SINGAPORE
\newline
Universit\'{e} Paris 7 - Math\'{e}matiques - Case 7012\quad\\
175/179 rue du Chevaleret\\
75013 Paris, FRANCE}
\email{berrick@math.nus.edu.sg
\newline
karoubi@math.jussieu.fr}
\begin{abstract}Rognes and Weibel used Voevodsky's work on the Milnor conjecture to deduce the
strong Dwyer-Friedlander form of the Lichtenbaum-Quillen conjecture at the
prime $2$. In consequence (the $2$-completion of) the classifying space for
algebraic $K$-theory of the integers $\mathbb{Z}[1/2]$ can be expressed as a
fiber product of well-understood spaces $BO$ and $B\mathrm{GL}(\mathbb{F}%
_{3})^{+}$ over $BU$. Similar results are now obtained for Hermitian
$K$-theory and the classifying spaces of the integral symplectic and
orthogonal groups. For the integers $\mathbb{Z}[1/2]$, this leads to
computations of the $2$-primary Hermitian $K$-groups and affirmation of the
Lichtenbaum-Quillen conjecture in the framework of Hermitian $K$-theory.
\end{abstract}
\maketitle

%
%
%
%
%
%
%
%
%
%
%
%
%
%
%
%
%
%
%
%
%
%
%
%
%
%
%
%
%
%
%
%
%
%
%
%
%
%
%
%
%
%
%
%
%
%
%
%
%
%
%
%
%
%
%
%
%
%
%
%
%
%
%
%
%
%
%
%
%
%
%
%
%

%
%
%
%
%
%
%
%
%
%
%
%
%
%
%
%
%
%
%
%
%
%
%
%
%
%
%
%
%
%
%
%
%
%
%
%
%
%
%
%
%
%
%
%
%
%
%
%
%
%
%
%
%
%
%
%
%
%
%
%
%
%
%
%
%
%
%
%
%
%
%
%
%

%
%
%
%
%
%
%
%
%
%
%
%
%
%
%
%
%
%
%
%
%
%
%
%
%
%
%
%
%
%
%
%
%
%
%
%
%
%
%
%
%
%
%
%
%
%
%
%
%
%
%
%
%
%
%
%
%
%
%
%
%
%
%
%
%
%
%
%
%
%
%
%
%

%
%
%
%
%
%
%
%
%
%
%
%
%
%
%
%
%
%
%
%
%
%
%
%
%
%
%
%
%
%
%
%
%
%
%
%
%
%
%
%
%
%
%
%
%
%
%
%
%
%
%
%
%
%
%
%
%
%
%
%
%
%
%
%
%
%
%
%
%
%
%
%
%

%
%
%
%
%
%
%
%
%
%
%
%
%
%
%
%
%
%
%
%
%
%
%
%
%
%
%
%
%
%
%
%
%
%
%
%
%
%
%
%
%
%
%
%
%
%
%
%
%
%
%
%
%
%
%
%
%
%
%
%
%
%
%
%
%
%
%
%
%
%
%
%
%

%
%
%
%
%
%
%
%
%
%
%
%
%
%
%
%
%
%
%
%
%
%
%
%
%
%
%
%
%
%
%
%
%
%
%
%
%
%
%
%
%
%
%
%
%
%
%
%
%
%
%
%
%
%
%
%
%
%
%
%
%
%
%
%
%
%
%
%
%
%
%
%
%

%
%
%
%
%
%
%
%
%
%
%
%
%
%
%
%
%
%
%
%
%
%
%
%
%
%
%
%
%
%
%
%
%
%
%
%
%
%
%
%
%
%
%
%
%
%
%
%
%
%
%
%
%
%
%
%
%
%
%
%
%
%
%
%
%
%
%
%
%
%
%
%
%

%
%
%
%
%
%
%
%
%
%
%
%
%
%
%
%
%
%
%
%
%
%
%
%
%
%
%
%
%
%
%
%
%
%
%
%
%
%
%
%
%
%
%
%
%
%
%
%
%
%
%
%
%
%
%
%
%
%
%
%
%
%
%
%
%
%
%
%
%
%
%
%
%

\pagestyle{myheadings}
\setcounter{section}{-1}%

\section{Introduction}

In \cite{Bok}, B\"{o}kstedt introduced the study of the commuting square%
\begin{equation}%
\begin{array}
[c]{ccc}%
B\mathrm{GL}(\mathbb{Z}^{\prime})_{\#}^{+} & \longrightarrow &  B\mathrm{GL}%
(\mathbb{R})_{\#}\\
\downarrow &  & \downarrow\\
B\mathrm{GL}(\mathbb{F}_{3})_{\#}^{+} & \overset{b}{\longrightarrow} &
B\mathrm{GL}(\mathbb{C})_{\#}%
\end{array}
\label{Diag 1}%
\end{equation}
Here $\mathbb{Z}^{\prime}$ denotes the ring $\mathbb{Z}[1/2]$, $\mathbb{F}%
_{3}$ the finite field with three elements, and, for $F=\mathbb{R},\mathbb{C}%
$, $B\mathrm{GL}(F\mathbb{)}$ is the classifying space of the infinite general
linear group $\mathrm{GL}(F)$ with the usual topology. The symbol $\#$
indicates the $2$-adic completion, and the map $b$ the Brauer lift,
corresponding to the fibring of Adams' map $\psi^{3}-1$ on $BU=B\mathrm{GL}%
(\mathbb{C})$. The remaining maps are induced from the obvious ring homomorphisms.

The Dwyer-Friedlander formulation of the Lichtenbaum-Quillen conjecture for
$\mathbb{Z}$ at the prime $2$ is that the above square is homotopy cartesian
\cite{DF}Conjecture 1.3, Proposition 4.2. This has been affirmed in work of
Rognes and Weibel \cite{RW} -- see \cite{W: CR}Corollary 8.

Since the homotopy fiber of the map $B\mathrm{GL}(\mathbb{R)}\rightarrow
B\mathrm{GL}(\mathbb{C)}$ is the homogeneous space $\mathrm{GL}(\mathbb{C)}%
$/$\mathrm{GL}(\mathbb{R)}$, which has the homotopy type of $\Omega
^{7}(B\mathrm{GL}(\mathbb{R}))$ by Bott periodicity \cite{Bott}, we may also
write the homotopy fibration%
\[
\Omega^{7}(B\mathrm{GL}(\mathbb{R))}_{\#}\longrightarrow B\mathrm{GL}%
(\mathbb{Z}^{\prime}\mathbb{)}_{\#}^{+}\longrightarrow B\mathrm{GL}%
(\mathbb{F}_{3}\mathbb{)}_{\#}^{+}\text{.}%
\]
The purpose of this paper is to prove analogous results (see Section 2 below)
with the orthogonal and symplectic groups over $\mathbb{Z}^{\prime}$
substituted for the general linear group. Our motivation comes from previous
work by the first author involving the mapping class group \cite{Be}, and by
the second author on the analogue of Bott periodicity in Hermitian $K$-theory
\cite{K:AnnM112 fun}. This subject is of course related to the computations by
A. Borel \cite{Bor} of the rational cohomology of arithmetic groups.

\medskip

\textit{Acknowledgements. }We thank A. Bak, L. Fajstrup, E.\thinspace M.
Friedlander, H. Hamraoui, J. Hornbostel, B. Kahn and L. N. Vaserstein for
their kind interest in this work. In this regard we would like to acknowledge
the conscientiousness of L. Fajstrup, whose recent amendment to some
calculations of \cite{Faj}\S8 provides independent confirmation of a key point
in our proof of Theorem C of Section \ref{main results}.

\section{\label{background}Motivational background}

The commutative diagram below appears in \cite{Be}. In it, $\mathrm{Br}_{g}$
denotes the $g$-strand braid group and the map \textrm{Artin }is Artin's
representation of $\mathrm{Br}_{g}$ as automorphisms of the free group
$\mathrm{Fr}_{g}$ on $g$ generators. The group $\Gamma_{g,1}$ is the mapping
class group of a surface of genus $g$ with one boundary component. The map
$\psi_{g}$ is constructed by Vershinin \cite{Vershinin}p.1000. $H$ is the
hyperbolic map sending a matrix $A$ to the matrix $\left(
\begin{array}
[c]{cc}%
A & O\\
O & ^{\mathrm{t}}A^{-1}%
\end{array}
\right)  $.
\[%
\begin{array}
[c]{ccccccc}%
&  & \Gamma_{g,1} & \overset{\mathrm{id}}{\longrightarrow} & \Gamma_{g,1} &  &
\\
& \ \quad\nearrow\psi_{g} &  &  &  & \ \ \searrow & \\
\mathrm{Br}_{g} & \ \ \twoheadrightarrow & \Sigma_{g} & \rightarrowtail &
O_{g}(\mathbb{Z}) & \overset{H}{\longrightarrow} & \mathrm{Sp}_{2g}%
(\mathbb{Z})\\
& \!_{_{\mathrm{Artin}}}\searrow\quad &  &  & \downarrow &  & \downarrow\\
&  & \mathrm{Aut}(\mathrm{Fr}_{g}) & \twoheadrightarrow & \mathrm{GL}%
_{g}(\mathbb{Z}) & \overset{H}{\longrightarrow} & \mathrm{GL}_{2g}(\mathbb{Z})
\end{array}
\]
Here all maps are injective except for the surjections $\mathrm{Aut}%
(\mathrm{Fr}_{g})\twoheadrightarrow\mathrm{GL}_{g}(\mathbb{Z})$,
$\mathrm{Br}_{g}\twoheadrightarrow\Sigma_{g}$ and $\Gamma_{g,1}%
\twoheadrightarrow\mathrm{Sp}_{2g}(\mathbb{Z})$.

Now, as in \cite{Be}, combine with the inclusions in the real symplectic and
general linear groups, stabilize, and take $B(\,\cdot\,)^{+}$, to get{\tiny
\[%
\begin{array}
[c]{ccccccccc}%
&  & B\Gamma_{\infty}^{+} & \overset{\mathrm{id}}{\rightarrow} &
B\Gamma_{\infty}^{+} &  &  &  & BU\medskip\\
& \nearrow &  &  &  & \searrow &  &  & \downarrow^{\simeq}\medskip\\
B\mathrm{Br}_{\infty}^{+} & \rightarrow &  B\Sigma_{\infty}^{+} & \rightarrow
&  BO(\mathbb{Z})^{+} & \rightarrow &  B\mathrm{Sp}(\mathbb{Z})^{+} &
\rightarrow &  B\mathrm{Sp}\mathbb{R}\\
& \searrow &  &  & \downarrow &  & \downarrow &  & \downarrow\\
&  &  B\mathrm{Aut}(\mathrm{Fr}_{\infty})^{+} & \rightarrow &  B\mathrm{GL}%
(\mathbb{Z})^{+} & \rightarrow &  B\mathrm{GL}(\mathbb{Z})^{+} & \rightarrow
&  B\mathrm{GL}\mathbb{R}%
\end{array}
\]
}The challenge is to describe these maps in homotopy theory. For instance,
\cite{CharneyCohen} discusses $B\Gamma_{\infty}^{+}\rightarrow B\mathrm{GL}%
(\mathbb{Z})^{+}$, while \cite{MadsenTillmann} looks at $B\Gamma_{\infty}%
^{+}\rightarrow BU$.

In this work we focus on $BO(\mathbb{Z}^{\prime})^{+}$ and $B\mathrm{Sp}%
(\mathbb{Z}^{\prime})^{+}$ as accessible approximations to $BO(\mathbb{Z}%
)^{+}$ and $B\mathrm{Sp}(\mathbb{Z})^{+}$. Indeed, if we extrapolate the
results known in the higher Hermitian $K$-theory of rings $A$ (where $2$ is
invertible) and the results known in lower $K$-theory, then it is reasonable
to conjecture that, up to odd finite torsion, for $i>1$ we have the
isomorphisms
\[
\pi_{i}(BO(\mathbb{Z})^{+})\cong\pi_{i}(BO(\mathbb{Z}^{\prime})^{+}%
)\quad\text{and\quad}\pi_{i}(B\mathrm{Sp}(\mathbb{Z})^{+})\cong\pi
_{i}(B\mathrm{Sp}(\mathbb{Z}^{\prime})^{+})\text{.}%
\]

\section{\label{main results}Main results}

We recall first some notations from \cite{K:AnnM112 fun}. Let $A$ be a ring
provided with an antiinvolution $x\mapsto\check{x}$, and let $\varepsilon$ be
an element of the center of $A$ such that $\varepsilon\check{\varepsilon}=1$.
We assume also the existence of an element $\lambda$ of the center such that
$\lambda+\check{\lambda}=1$. In most cases, $\lambda$ $=1/2$ and
$\varepsilon=\pm1$. In this setting, a central role is played by the
$\varepsilon$\emph{-orthogonal group}, which is the group of automorphisms of
the $\varepsilon$-hyperbolic module ${}_{\varepsilon}H(A^{n})$, denoted by
${}_{\varepsilon}O_{n,n}(A)$: its elements can be described as $2\times2$
matrices written in $n$-blocks%
\[
M=\left(
\begin{array}
[c]{cc}%
a & b\\
c & d
\end{array}
\right)
\]
such that $M^{\ast}M=MM^{\ast}=I$, where the `$\varepsilon$-hyperbolic
adjoint' $M^{\ast}$ is defined as
\[
M^{\ast}=\left(
\begin{array}
[c]{cc}%
^{\mathrm{t}}\check{d} & \varepsilon\,^{\mathrm{t}}\check{b}\\
\check{\varepsilon}\,^{\mathrm{t}}\check{c} & ^{\mathrm{t}}\check{a}%
\end{array}
\right)
\]
given by conjugating the involute transpose of $M$ by the matrix%
\[
{}_{\varepsilon}J_{n}=\left(
\begin{array}
[c]{cc}%
0 & \varepsilon I_{n}\\
I_{n} & 0
\end{array}
\right)  \text{.}%
\]
For instance, when $A$ is commutative (corresponding to ${}\check{\cdot
}\,=\mathrm{id}$ also being an involution), the usual symplectic group\textrm{
}$\mathrm{Sp}_{2n}(A)$ is just ${}_{-1}O_{n,n}(A)$ with our notations. As
other examples, the classical orthogonal groups ${}_{1}O_{n,n}(\mathbb{R)}$
and ${}_{1}O_{n,n}(\mathbb{C)}$ have the homotopy type of $O(n)\times O(n)$
and $O(2n)$ respectively, whereas $\mathrm{Sp}_{2n}(\mathbb{R)}$ has the
homotopy type of the unitary group $U(n)$ (see, for example \cite{Helgason}).
We denote by ${}_{\varepsilon}O(A)$ the direct limit of the ${}_{\varepsilon
}O_{n,n}(A)$ with respect to the obvious inclusions within each of the four
component blocks. If $A=\mathbb{R}$ or $\mathbb{C}$, we provide ${}%
_{\varepsilon}O(A)$ with its usual topology unless otherwise stated.

The first aim of this paper is to prove the following squares homotopy
cartesian:%
\begin{equation}%
\begin{array}
[c]{ccc}%
B_{1}O(\mathbb{Z}^{\prime})_{\#}^{+} & \longrightarrow &  B_{1}O(\mathbb{R}%
)_{\#}\\
\downarrow &  & \downarrow\\
B_{1}O(\mathbb{F}_{3})_{\#}^{+} & \longrightarrow &  B_{1}O(\mathbb{C})_{\#}%
\end{array}
\label{Diag 2}%
\end{equation}
and%
\begin{equation}%
\begin{array}
[c]{ccc}%
B_{-1}O(\mathbb{Z}^{\prime})_{\#}^{+} & \longrightarrow &  B_{-1}%
O(\mathbb{R})_{\#}\\
\downarrow &  & \downarrow\\
B_{-1}O(\mathbb{F}_{3})_{\#}^{+} & \longrightarrow &  B_{-1}O(\mathbb{C})_{\#}%
\end{array}
\label{Diag 3}%
\end{equation}
Since we have the well-known homotopy equivalences
\[
_{1}O(\mathbb{R})\simeq O\times O,\qquad_{1}O(\mathbb{C})\simeq O,\qquad
_{-1}O(\mathbb{R})\simeq U,\qquad_{-1}O(\mathbb{C})\simeq\mathrm{Sp},
\]
on the one hand, and (cf. \cite{K:LNM343}p.311)%
\[
\mathrm{Sp}/U\simeq\Omega^{6}(BO)
\]
on the other, we obtain homotopy fibrations%
\[%
\begin{array}
[c]{ccccc}%
BO_{\#} & \longrightarrow &  B_{1}O(\mathbb{Z}^{\prime})_{\#}^{+} &
\longrightarrow &  B_{1}O(\mathbb{F}_{3})_{\#}^{+}%
\end{array}
\]%
\[%
\begin{array}
[c]{ccccc}%
\Omega^{6}(BO)_{\#} & \longrightarrow &  B\mathrm{Sp}(\mathbb{Z}^{\prime
})_{\#}^{+} & \longrightarrow &  B\mathrm{Sp}(\mathbb{F}_{3})_{\#}^{+}%
\end{array}
\]

Before starting to prove these statements, we should put them in a slightly
more general framework. We have to consider the full classifying spaces of
algebraic $K$-theory and Hermitian $K$-theory which we denote by
$\mathcal{K}(A)$ and $_{\varepsilon}\mathcal{L}(A)$ respectively, following
the notations of \cite{K:AnnM112 fun}\footnote{These notations are different
from the ones generally used in surgery theory.}. Specifically, we have
homotopy equivalences%
\[
\mathcal{K}(A)\simeq K_{0}(A)\times B\mathrm{GL}(A)^{+}\qquad\text{and\qquad
}{}_{\varepsilon}\mathcal{L}(A)\simeq{}_{\varepsilon}L_{0}(A)\times
B{}_{\varepsilon}O(A)^{+}.
\]
Here $K_{0}(A)$ denotes the usual Grothendieck group of the category of
finitely generated projective $A$-modules and ${}_{\varepsilon}L_{0}(A)$ the
`Witt-\-Grothen\-dieck group' of the category constructed from the same
objects with the extra structure of a nondegenerate $\varepsilon$-Hermitian
form. There are obvious modifications when $A=\mathbb{R}$ or $\mathbb{C}$. A
more functorial encoding of the $\pi_{0}$ information, which we need for
consideration of fixed points under the involution, is given by the
formulations $\mathcal{K}(A)=\Omega(B\mathrm{GL}(SA)^{+})$ and $_{\varepsilon
}\mathcal{L}(A)=\Omega(B_{\varepsilon}O(SA)^{+})$ -- see Appendix A. For our
notation involving the $2$-adic completion $X_{\#}$ of a non-connected space
$X$, we use the convention that when, as here, $X$ is the loop space $X=\Omega
Y$ of a connected space $Y$, then $X_{\#}$ is the loop space $\Omega(Y_{\#})$.

For instance, the first diagram (\ref{Diag 1}) may be written in an equivalent
way as
\begin{equation}%
\begin{array}
[c]{ccc}%
\mathcal{K}(\mathbb{Z}^{\prime})_{\#} & \longrightarrow & \mathcal{K}%
(\mathbb{R})_{\#}\\
\downarrow &  & \downarrow\\
\mathcal{K}(\mathbb{F}_{3})_{\#} & \longrightarrow & \mathcal{K}%
(\mathbb{C})_{\#}%
\end{array}
\label{Diag 1'}%
\end{equation}
(we just cross the spaces by $K_{0}(A)_{\#}=\mathbb{Z}_{\#}$). The second and
third diagrams may likewise be written in the following form (with
$\varepsilon=\pm1$):%
\begin{equation}%
\begin{array}
[c]{ccc}%
_{\varepsilon}\mathcal{L}(\mathbb{Z}^{\prime})_{\#} & \longrightarrow &
_{\varepsilon}\mathcal{L}(\mathbb{R})_{\#}\\
\downarrow &  & \downarrow\\
_{\varepsilon}\mathcal{L}(\mathbb{F}_{3})_{\#} & \longrightarrow &
_{\varepsilon}\mathcal{L}(\mathbb{C})_{\#}%
\end{array}
\label{Diag 2/3}%
\end{equation}
This is clear for $\varepsilon=-1$ since $_{-1}L_{0}(A)\cong\mathbb{Z}$ for
$A=\mathbb{Z}^{\prime}$, $\mathbb{F}_{3}$, $\mathbb{R}$ or $\mathbb{C}$, as
detailed in (\ref{-1L0}) below. For $\varepsilon=+1$ it requires a little more
care, and uses the fact, proven in (\ref{1L0}) below, that the following
square of Witt-Grothendieck groups is cartesian:%
\[%
\begin{array}
[c]{ccc}%
_{1}L_{0}(\mathbb{Z}^{\prime})\cong\mathbb{Z\oplus Z\oplus Z}/2 &
\longrightarrow & \mathbb{Z\oplus Z}\cong{}_{1}L_{0}(\mathbb{R})\\
\downarrow &  & \downarrow\\
_{1}L_{0}(\mathbb{F}_{3})\cong\mathbb{Z\oplus Z}/2 & \longrightarrow &
\mathbb{Z}\cong{}_{1}L_{0}(\mathbb{C})
\end{array}
\]
In this setting, we can now state our first main theorem.

\smallskip

\noindent\textbf{Theorem A. }\label{L cartesian}\emph{Diagram \ref{Diag 2/3}
above (for }$\varepsilon=\pm1$\emph{)}
\[%
\begin{array}
[c]{ccc}%
_{\varepsilon}\mathcal{L}(\mathbb{Z}^{\prime})_{\#} & \longrightarrow &
_{\varepsilon}\mathcal{L}(\mathbb{R})_{\#}\\
\downarrow &  & \downarrow\\
_{\varepsilon}\mathcal{L}(\mathbb{F}_{3})_{\#} & \longrightarrow &
_{\varepsilon}\mathcal{L}(\mathbb{C})_{\#}%
\end{array}
\]
\emph{is homotopy cartesian.}

\emph{\smallskip}

In consequence, we calculate the groups ${}_{\varepsilon}L_{i}(\mathbb{Z}%
^{\prime})$, up to finite odd order subgroups. For comparison we also include
$K_{i}(\mathbb{Z}^{\prime})$ information, which follows immediately from what
is known about $K_{i}(\mathbb{Z})$ \cite{W: CR}, and the localization exact
sequence%
\[
K_{i+1}(\mathbb{Z})\rightarrow K_{i+1}(\mathbb{Z}^{\prime})\rightarrow
K_{i}(\mathbb{F}_{2})\rightarrow K_{i}(\mathbb{Z})\rightarrow K_{i}%
(\mathbb{Z}^{\prime})\text{,}%
\]
given that $K_{i}(\mathbb{F}_{2})$ is a finite group of odd order for $i>0$.

\smallskip

\smallskip

\noindent\textbf{Theorem B. }\emph{Modulo a finite group of odd order, the
groups }$K_{i}(\mathbb{Z}^{\prime})$\emph{ and }${}_{\varepsilon}%
L_{i}(\mathbb{Z}^{\prime})$\emph{ for }$i\geq0$\emph{ are as follows, where
}$\delta_{i0}$\emph{ denotes the Kronecker delta, and }$2^{t}$\emph{ is the
}$2$\emph{-primary part of }$i+1$\emph{.}
\[%
\begin{array}
[c]{cccc}%
i\;(\mathrm{\operatorname{mod}\,}8) & K_{i}(\mathbb{Z}^{\prime}) & _{-1}%
L_{i}(\mathbb{Z}^{\prime}) & _{1}L_{i}(\mathbb{Z}^{\prime})\\
0 & \delta_{i0}\mathbb{Z} & \delta_{i0}\mathbb{Z} & \delta_{i0}\mathbb{Z\oplus
Z\oplus Z}/2\\
1 & \mathbb{Z\oplus Z}/2 & 0 & \mathbb{Z}/2\oplus\mathbb{Z}/2\oplus
\mathbb{Z}/2\\
2 & \mathbb{Z}/2 & \mathbb{Z} & \mathbb{Z}/2\oplus\mathbb{Z}/2\\
3 & \mathbb{Z}/16 & \mathbb{Z}/16 & \mathbb{Z}/8\\
4 & 0 & \mathbb{Z}/2 & \mathbb{Z}\\
5 & \mathbb{Z} & \mathbb{Z}/2 & 0\\
6 & 0 & \mathbb{Z} & 0\\
7 & \mathbb{Z}/2^{t+1} & \mathbb{Z}/2^{t+1} & \mathbb{Z}/2^{t+1}%
\end{array}
\]

A further consequence is affirmation of the `Lichtenbaum-Quillen conjecture'
for Hermitian $K$-theory, as follows. (Here the notation ${}_{\varepsilon
}\mathbb{Z}/2$ keeps record of the particular $\mathbb{Z}/2$ action. Thus the
homotopy fixed point set
\[
\mathcal{K}(\mathbb{Z}^{\prime}){}^{h({}_{\varepsilon}\mathbb{Z}%
/2)}:=\mathrm{map}_{{}_{\varepsilon}\mathbb{Z}/2}(E\mathbb{Z}/2,\,\mathcal{K}%
(\mathbb{Z}^{\prime}))
\]
is the space of maps equivariant under the ${}_{\varepsilon}\mathbb{Z}/2$
action, where $E\mathbb{Z}/2$ is a contractible free $\mathbb{Z}/2$-space,
such as $S^{\infty}$ with antipodal action, and $\mathbb{Z}/2$ acts on
$\mathcal{K}(\mathbb{Z}^{\prime})$ via conjugation by the matrices
${}_{\varepsilon}J_{n}$ as described above.)

\smallskip

\smallskip

\noindent\textbf{Theorem C. }\emph{For} $\varepsilon=\pm1$\emph{, the natural
map}%
\[
{}_{\varepsilon}\mathcal{L}(\mathbb{Z}^{\prime})=\mathcal{K}(\mathbb{Z}%
^{\prime})^{{}_{\varepsilon}\mathbb{Z}/2}\longrightarrow\mathcal{K}%
(\mathbb{Z}^{\prime})^{h({}_{\varepsilon}\mathbb{Z}/2)}\text{,}%
\]
\emph{from the fixed point set to the homotopy fixed point set of the}
${}_{\varepsilon}\mathbb{Z}/2$ \emph{action on }$\mathcal{K}(\mathbb{Z}%
^{\prime})$\emph{, becomes a homotopy equivalence after }$2$\emph{-adic completion.}

For comparison, we remark that in \cite{Rognes} it is shown that, for the
spectrum $T(\mathbb{Z})$ for topological Hochschild homology of the integers,
the natural map from the fixed point set to the homotopy fixed point set,
$T(\mathbb{Z})^{\mathbb{Z}/2}\rightarrow T(\mathbb{Z})^{h\mathbb{Z}/2}$, also
induces a homotopy equivalence of $2$-completed connective covers.

Our results also hold at the level of spectra. For, as in \cite{Rognes JPAA
1999}, Diagram \ref{Diag 1} above may be recast as a diagram of $2$-completed
connective spectra. Then, because in order to achieve the passage from
algebraic to Hermitian $K$-theory we arrange that maps of infinite loop spaces
are induced by maps of rings and so are infinitely deloopable, our results
similarly can be expressed in terms of $2$-completed connective spectra. (We
note that it is essential to consider connective spectra, since the negative
$K$-groups of any discrete regular Noetherian ring like $\mathbb{Z}^{\prime}$
and $\mathbb{F}_{3}$ are all zero, according to a well-known theorem of Bass;
in contrast, the negative topological $K$-groups of $\mathbb{R}$ and
$\mathbb{C}$, by periodicity, do not always coincide. By our induction results
below, similar facts hold for Hermitian $K$-groups in negative dimensions.)\smallskip

\textit{Strategy of the proofs, and organization of the paper. } Let us denote
by $_{\varepsilon}\overline{\mathcal{L}}(\mathbb{Z}^{\prime})$ the homotopy
cartesian product of $_{\varepsilon}\mathcal{L}(\mathbb{F}_{3})$ and
$_{\varepsilon}\mathcal{L}(\mathbb{R})$ over $_{\varepsilon}\mathcal{L}%
(\mathbb{C})$. Tautologically, for Theorem A we want to prove that the map%
\[
{}_{\varepsilon}\phi:\,_{\varepsilon}\mathcal{L(}\mathbb{Z}^{\prime}%
)_{\#}\longrightarrow\,_{\varepsilon}\overline{\mathcal{L}}(\mathbb{Z}%
^{\prime})_{\#}%
\]
is a homotopy equivalence. Our strategy of proof is as follows:

The first step is to check that ${}_{\varepsilon}\phi$ induces an isomorphism
on $\pi_{0}$ and $\pi_{1}$ for $\varepsilon=\pm1$. This requires verifications
that are detailed in Section \ref{low} below.

As discussed at the outset, we have the homotopy cartesian diagram (\ref{Diag
1'}):%
\[%
\begin{array}
[c]{ccc}%
\mathcal{K}(\mathbb{Z}^{\prime})_{\#} & \longrightarrow & \mathcal{K}%
(\mathbb{R})_{\#}\\
\downarrow &  & \downarrow\\
\mathcal{K}(\mathbb{F}_{3})_{\#} & \longrightarrow & \mathcal{K}%
(\mathbb{C})_{\#}%
\end{array}
\]
In other words, if we denote by $\overline{\mathcal{K}}(\mathbb{Z}^{\prime})$
the homotopy cartesian product of $\mathcal{K}(\mathbb{F}_{3})$ and
$\mathcal{K}(\mathbb{R})$ over $\mathcal{K}(\mathbb{C})$, the obvious map%
\[
\psi:\mathcal{K}(\mathbb{Z}^{\prime})_{\#}\longrightarrow\overline
{\mathcal{K}}(\mathbb{Z}^{\prime})_{\#}%
\]
is a homotopy equivalence.

The last step of the proof, in Section \ref{htpy equiv Thm A}, is to show that
the verifications above and the fact that $\psi$ is a homotopy equivalence
imply that ${}_{\varepsilon}\phi$ is also a homotopy equivalence. This step is
similar in spirit to the argument used already in \cite{K2}, and is discussed
in Section \ref{induction methods}, together with several other `Karoubi
induction' methods that we construct for later use. Here however it requires a
careful study of the map $\mathcal{K}(\mathbb{F}_{3})\longrightarrow
\mathcal{K}(\mathbb{C})$ (the Brauer lifting) and its analogue in Hermitian
$K$-theory, in order that required equivalences be induced from ring
homomorphisms. A further ingredient is that rational and mod\thinspace$2$
arguments need to be handled separately.

Theorem B then follows readily from Theorem A by consideration in Section
\ref{L groups} of the various exact sequences of homotopy groups that arise
from homotopy cartesian squares.

In order to prove Theorem C, we need to compare the homotopy cartesian square
(\ref{Diag 2/3}) of Theorem A with the homotopy fixed point sets that occur in
the homotopy cartesian square of Diagram \ref{Diag 1'}. Again, the behaviour
of $\mathcal{K}(\mathbb{F}_{3})$ requires some care. However, the most
treacherous part of the investigation concerns the equivalence ${}%
_{\varepsilon}\mathcal{L}(\mathbb{R})\longrightarrow\mathcal{K}(\mathbb{R}%
)^{h({}_{\varepsilon}\mathbb{Z}/2)}$. Our treatment draws on the
representation of topological $K$-theory by means of Fredholm operators on a
Hilbert space, and leads to a result of independent interest, on the algebraic
$K$-theory of simple algebras with involution.

\section{\label{induction methods}Induction methods}

\begin{theorem}
\label{K-L bootstrap}Let $n_{0}\in\mathbb{Z}$ be arbitrary. Suppose that
$\gamma:A\longrightarrow B$ is a Hermitian ring map of $\mathbb{Z}^{\prime}%
$-algebras inducing an isomorphism%
\[
{}_{\varepsilon}L_{n}(\gamma):{}_{\varepsilon}L_{n}(A)\overset{\cong
}{\longrightarrow}{}_{\varepsilon}L_{n}(B)
\]
for $\varepsilon=\pm1$ and $n=n_{0},\,n_{0}+1$.

\noindent\emph{(a) Upward induction.}\textrm{ }Let $r\geq1$, and suppose that
$\gamma$ induces an isomorphism%
\[
K_{n}(\gamma):K_{n}(A)\overset{\cong}{\longrightarrow}K_{n}(B)
\]
whenever $n_{0}\leq n\leq n_{0}+r$. Then $\gamma$ also induces an isomorphism
\[
_{\varepsilon}L_{n}(\gamma):{}_{\varepsilon}L_{n}(A)\overset{\cong
}{\longrightarrow}{}_{\varepsilon}L_{n}(B)
\]
whenever $\varepsilon=\pm1$ and $n_{0}\leq n\leq n_{0}+r$.

\noindent\emph{(b) Downward induction.}\textrm{ }Let $s\geq1$, and suppose
that $\gamma$ induces an isomorphism%
\[
K_{n}(\gamma):K_{n}(A)\overset{\cong}{\longrightarrow}K_{n}(B)
\]
whenever $n_{0}-s\leq n\leq n_{0}$ and an epimorphism when $n=n_{0}+1$. Then
$\gamma$ also induces an isomorphism
\[
_{\varepsilon}L_{n}(\gamma):{}_{\varepsilon}L_{n}(A)\overset{\cong
}{\longrightarrow}{}_{\varepsilon}L_{n}(B)
\]
whenever $\varepsilon=\pm1$ and $n_{0}-s\leq n\leq n_{0}$.
\end{theorem}

\noindent\textbf{Remark.\quad} It will be clear from the proof that the
assertion is also valid for homotopy groups with finite coefficients, and
modulo a Serre class of abelian groups.

\noindent\textbf{Proof.\quad} \textbf{(a)}\quad We prove the theorem by
induction on $r$, the case $r=1$ being given. Suppose therefore that the
theorem holds for a given $r\geq1$. Following the notations of \cite{K:AnnM112
fun}, we have for all $n$ a commutative diagram of exact sequences%
\[%
\begin{array}
[c]{ccccccccc}%
_{\varepsilon}L_{n+1}(A) & \longrightarrow &  K_{n+1}(A) & \longrightarrow &
_{\varepsilon}V_{n}(A) & \longrightarrow & _{\varepsilon}L_{n}(A) &
\longrightarrow &  K_{n}(A)\\
\, &  &  &  &  &  &  &  & \\
\downarrow^{_{\varepsilon}L_{n+1}(\gamma)} &  & \downarrow^{K_{n+1}(\gamma)} &
& \downarrow^{\varphi_{n}} &  & \downarrow^{_{\varepsilon}L_{n}(\gamma)} &  &
\downarrow^{K_{n}(\gamma)}\\
&  &  &  &  &  &  &  & \\
_{\varepsilon}L_{n+1}(B) & \longrightarrow &  K_{n+1}(B) & \longrightarrow &
_{\varepsilon}V_{n}(B) & \longrightarrow & _{\varepsilon}L_{n}(B) &
\longrightarrow &  K_{n}(B)
\end{array}
\]

Now write $m=n_{0}+r$. The five lemma and the induction hypothesis imply that
$\varphi_{m-1}$ is an isomorphism and $\varphi_{m}$ is an epimorphism. By the
fundamental theorem in Hermitian $K$-theory \cite{K:AnnM112 fun}, the groups
${}_{\varepsilon}V_{n}$ and ${}_{-\varepsilon}U_{n+1}$ are naturally
isomorphic. Therefore, we can write a second commutative diagram of exact
sequences:{\small
\[%
\begin{array}
[c]{ccccccccc}%
{}_{-\varepsilon}U_{m+1}(A) & \rightarrow &  K_{m+1}(A) & \rightarrow &
{}_{-\varepsilon}L_{m+1}(A) & \rightarrow & {}_{-\varepsilon}U_{m}(A) &
\rightarrow &  K_{m}(A)\\
&  &  &  &  &  &  &  & \\
\downarrow^{\varphi_{m}} &  & \downarrow^{K_{m+1}(\gamma)} &  & \downarrow
^{{}_{-\varepsilon}L_{m+1}(\gamma)} &  & \downarrow^{\varphi_{m-1}} &  &
\downarrow^{K_{m}(\gamma)}\\
&  &  &  &  &  &  &  & \\
{}_{-\varepsilon}U_{m+1}(B) & \rightarrow &  K_{m+1}(B) & \rightarrow &
{}_{-\varepsilon}L_{m+1}(B) & \rightarrow & {}_{-\varepsilon}U_{m}(B) &
\rightarrow &  K_{m}(B)
\end{array}
\]
}According to the above, $K_{m+1}(\gamma)$, $\varphi_{m-1}$, $K_{m}(\gamma)$
are isomorphisms and $\varphi_{m}$ is an epimorphism. Therefore,
${}_{-\varepsilon}L_{m+1}(\gamma)$ is an isomorphism by the five lemma.

\noindent\textbf{(b)}\quad Given that ${}_{\pm1}L_{n}(\gamma),{}_{\pm1}%
L_{n+1}(\gamma),K_{n-1}(\gamma)$ and $K_{n}(\gamma)$ are isomorphisms, and
that $K_{n+1}(\gamma)$ is an epimorphism, applying the five lemma to the map
of exact sequences (with vertical arrows and $\gamma$ omitted for brevity)%
\[%
\begin{array}
[c]{ccccccccccccc}%
K_{n+1} & \  & {}_{-\varepsilon}L_{n+1} & \  & \varphi_{n-1} & \  & K_{n} &
\  & {}_{-\varepsilon}L_{n} & \  & \varphi_{n-2} & \  & K_{n-1}%
\end{array}
\]
reveals that $\varphi_{n-1}\ $is an isomorphism and $\varphi_{n-2}$ a
monomorphism. Then the five lemma for the map of exact sequences%
\[%
\begin{array}
[c]{ccccccccc}%
K_{n} & \  & \varphi_{n-1} & \  & {}_{\varepsilon}L_{n-1} & \  & K_{n-1} &
\  & \varphi_{n-2}%
\end{array}
\]
gives ${}_{\varepsilon}L_{n-1}(\gamma)$ as an isomorphism, and so the
induction proceeds.\hfill$\Box\medskip$

We now present several consequences of this result that we shall need later,
and which are of independent interest. The first applications are to
mod$\,\ell$ $L$-theory, where $\ell\geq2$. To simplify notation, for
$n\in\mathbb{Z}$ we write
\[
{}_{\varepsilon}\tilde{L}_{n}(A)={}_{\varepsilon}L_{n}(A;\,\mathbb{Z}%
/\ell)\text{,}%
\]
being the $n\,$th $\mathbb{Z}/\ell$-homotopy group of the spectrum
$_{\varepsilon}\mathcal{L}(A)$ continued into negative dimensions by iterated
suspension; likewise for mod$\,\ell$ $K$-groups. The first result is a
technical device for moving from integral to mod$\,\ell$ $L$-theory.

\begin{Corollary}
\label{L to mod p L}In the notation of the theorem above, suppose that, for
both $\varepsilon=1$ and $\varepsilon=-1$, $\gamma$ induces isomorphisms on
${}_{\varepsilon}L_{n_{0}},$ ${}_{\varepsilon}L_{n_{0}+1},K_{n_{0}-1}$ and
$K_{n_{0}}$, and an epimorphism on $K_{n_{0}+1}$. If $\gamma$ induces
isomorphisms on the mod$\,\ell$ $K$-groups $\tilde{K}_{n}$ for $n_{0}-s\leq
n\leq n_{0}+r$, then it also induces isomorphisms on all $\mathrm{mod}\;\ell$
$L$-groups ${}_{\varepsilon}\tilde{L}_{n}$ in this range.
\end{Corollary}

\noindent\textbf{Proof.\quad} From (b) of the theorem, we obtain an
isomorphism on ${}_{\varepsilon}L_{n_{0}-1}$. Then the universal coefficients
formula%
\[
{}_{\varepsilon}L_{n}(A)\otimes\mathbb{Z}/\ell\rightarrowtail{}_{\varepsilon
}\tilde{L}_{n}(A)\twoheadrightarrow{}_{\ell}({}_{\varepsilon}L_{n-1}(A))
\]
with $n=n_{0},n_{0}+1$ gives the isomorphism on ${}_{\varepsilon}\tilde
{L}_{n_{0}},{}_{\varepsilon}\tilde{L}_{n_{0}+1}$ that enables both upward and
downward induction to commence.\hfill$\Box\medskip$

Next, we have a Mayer-Vietoris sequence.

\begin{Theorem}
\label{M-V}For a prime $\ell\geq2$, let
\[%
\begin{array}
[c]{ccc}%
A_{1} & \overset{\alpha}{\twoheadrightarrow} & A_{2}\\
\downarrow & ^{\cdot}\lrcorner & \downarrow\\
A_{3} & \overset{\beta}{\twoheadrightarrow} & A_{4}%
\end{array}
\]
be a cartesian square of Hermitian rings and homomorphisms, with $\beta$ (and
therefore $\alpha$) surjective, where each ring is a $\mathbb{Z}[1/2\ell
]$-algebra. Then for $\varepsilon\in\{\pm1\}$and all integers $n$ there is a
natural exact sequence
\[
\cdots\rightarrow{}_{\varepsilon}\tilde{L}_{n+1}(A_{4})\rightarrow
{}_{\varepsilon}\tilde{L}_{n}(A_{1})\rightarrow{}_{\varepsilon}\tilde{L}%
_{n}(A_{2})\oplus{}_{\varepsilon}\tilde{L}_{n}(A_{3})\rightarrow
{}_{\varepsilon}\tilde{L}_{n}(A_{4})\rightarrow\cdots
\]
\end{Theorem}

\noindent\textbf{Proof.\quad} As usual, the Mayer-Vietoris sequence follows
from an excision isomorphism between the relative terms obtained from the
surjections $\alpha$ and $\beta$. To obtain this isomorphism, first recall
from \cite{K:AnnM112 fun}p.262 that associated to a ring map $\varphi
:B\rightarrow D$ there is the fiber product $\Gamma(\varphi)$%
\[%
\begin{array}
[c]{ccc}%
\Gamma(\varphi) & \overset{}{\twoheadrightarrow} & CD\\
\downarrow & ^{\cdot}\lrcorner & \downarrow\\
SB & \overset{\varphi}{\twoheadrightarrow} & SD
\end{array}
\]
where $C$ and $S$ stand for cone and suspension respectively. This ring has
the property that $\Omega\mathcal{K}(\Gamma(\varphi))$ and $\Omega
{}_{\varepsilon}\mathcal{L(}\Gamma(\varphi))$ have the homotopy type of the
homotopy fibers of $\mathcal{K}(\varphi)$ and ${}_{\varepsilon}\mathcal{L(}%
\varphi)$. So here it suffices to check that the map $\Gamma(\alpha
)\rightarrow\Gamma(\beta)$ induces an isomorphism of ${}_{\varepsilon}%
\tilde{L}_{\ast}$-groups. This is an application of the preceding theorem to
the data that the excision isomorphism for $\mathbb{Z}[1/\ell]$-algebras
holds: first, for mod\thinspace$\ell$ $K$-theory (using surjectivity) by
\cite{W}; and, second, for negatively indexed ${}_{\varepsilon}\tilde{L}%
_{\ast}$-groups by application of universal coefficients to the results of
\cite{KV2}. \hfill$\Box\medskip$

Homotopy invariance also follows.

\begin{Theorem}
\label{htpy invar}For a prime $\ell\geq2$, a $\mathbb{Z}[1/2\ell]$-algebra $A$
and $\varepsilon\in\{\pm1\}$, there is a natural isomorphism%
\[
{}_{\varepsilon}\tilde{L}_{\ast}(A[x])\cong{}_{\varepsilon}\tilde{L}_{\ast
}(A)\text{.}%
\]
\end{Theorem}

\noindent\textbf{Proof.\quad} The corresponding result for mod\thinspace$\ell$
$K$-theory is proved in \cite{W}. Again, we check the result for negatively
indexed $_{\varepsilon}\tilde{L}_{\ast}$-groups. Since for suspension rings we
have from Appendix A that $S(A[x])\cong(SA)[x]$, the result is true in all
negative dimensions provided it holds for ${}_{\varepsilon}\tilde{L}_{0}$. As
the injection $A\rightarrowtail A[x]$ is split, the universal coefficients
formula%
\[
K_{0}(A)\otimes\mathbb{Z}/\ell\rightarrowtail\tilde{K}_{0}%
(A)\twoheadrightarrow{}_{\ell}K_{-1}(A)
\]
shows that homotopy invariance holds for $\tilde{K}_{0}$ precisely when it
holds for both $\mathbb{Z}/\ell\otimes K_{0}$ and ${}_{\ell}K_{-1}$; likewise
for ${}_{\varepsilon}\tilde{L}_{0}$. From homotopy invariance of
${}_{\varepsilon}W_{0}=\mathrm{Coker}(K_{0}\rightarrow{}_{\varepsilon}L_{0})$
established in \cite{K: AnnScENS1975} (and so, by suspension, of
${}_{\varepsilon}W_{-1}$), we use the invariance of $\mathbb{Z}/\ell\otimes
K_{0}$ (respectively ${}_{\ell}K_{-1}$) to deduce the invariance of
$\mathbb{Z}/\ell\otimes{}_{\varepsilon}L_{0}$ (respectively ${}_{\ell}%
({}_{\varepsilon}L_{-1})$). Hence we have homotopy invariance for
${}_{\varepsilon}\tilde{L}_{0}$ too. \hfill$\Box\medskip$

Here is another useful application of the methodology.

\begin{Proposition}
Let $A$ be a ring with $2$ invertible in $A$, and let $r\geq1$, $s\geq0$.
Suppose that the groups $K_{n}(A)$ (whenever $n_{0}-s\leq n\leq n_{0}+r$) and
${}_{\varepsilon}L_{i}(A)$ (for $\varepsilon\in\{\pm1\},$ $n\in\{n_{0}%
,n_{0}+1\}$) are finitely generated. Then the groups ${}_{\varepsilon}%
L_{n}(A)$ are finitely generated whenever $n_{0}-s\leq n\leq n_{0}+r$.
\end{Proposition}

\noindent\textbf{Proof.\quad} The result follows from Theorem \ref{K-L
bootstrap} above by considering the embedding of $A$ in its cone ring $CA$,
which has trivial $K$- and $L$-theory. (Clearly the argument works more
generally with finite generation replaced by any Serre class of abelian
groups.)\hfill$\Box\medskip$

This of course notably applies to Quillen's result \cite{Q: fg} on finite
generation of $K$-groups of rings of $S$-integers.

\begin{Corollary}
\label{fg}Let $A$ be a ring of $S$-integers in a number field such that
$\mathbb{Z}^{\prime}\subseteq A$. Then for all integers $n$ the groups
${}_{\varepsilon}L_{n}(A)$ are finitely generated.\hfill$\Box$
\end{Corollary}

\section{\label{low}Low-dimensional checks}

We write ${}_{\varepsilon}L_{n}(A)$ (resp. $_{\varepsilon}L_{n}^{\mathrm{top}%
}(A)$) for the homotopy groups $\pi_{n}({}_{\varepsilon}\mathcal{L}(A))$ in
the discrete case (resp. in the continuous case). Analogously, we write
$K_{n}(A)$ (resp. $K_{n}^{\mathrm{top}}(A)$) for the homotopy groups $\pi
_{n}(\mathcal{K}(A))$ in the same situations. Finally, $_{\varepsilon
}\overline{L}_{n}(\mathbb{Z}^{\prime})$ denotes the homotopy groups $\pi
_{n}(_{\varepsilon}\overline{\mathcal{L}}(\mathbb{Z}^{\prime}))$.

\begin{Void}
\label{-1L0}\textbf{The map} ${}_{-1}\phi:{}_{-1}\mathcal{L}(\mathbb{Z}%
^{\prime})_{\#}\longrightarrow{}_{-1}\overline{\mathcal{L}}(\mathbb{Z}%
^{\prime})_{\#}$ \textbf{induces an isomorphism on} $\pi_{0}$.
\end{Void}

In fact, this result holds even before $2$-completion. For, because all four
rings $A=\mathbb{Z}^{\prime}$, $\mathbb{F}_{3}$, $\mathbb{R}$, $\mathbb{C}$
are Dedekind, each $\pi_{0}(_{-1}\mathcal{L}(A))$ equals $\mathbb{Z}$,
detected by the (even) rank of the free symplectic $A$-inner product space
\cite{MH}p.7. Since each map of rings preserves rank, there is a cartesian
square%
\[%
\begin{array}
[c]{ccc}%
_{-1}L_{0}(\mathbb{Z}^{\prime})\cong\mathbb{Z} & \overset{\cong}%
{\longrightarrow} & \mathbb{Z}\cong{}_{-1}L_{0}(\mathbb{R})\\
\downarrow^{\cong} &  & \downarrow^{\cong}\\
_{-1}L_{0}(\mathbb{F}_{3})\cong\mathbb{Z} & \overset{\cong}{\longrightarrow} &
\mathbb{Z}\cong{}_{-1}L_{0}(\mathbb{C})
\end{array}
\]

On the other hand, the exact sequence that enables us to compute $\pi
_{0}(_{-1}\overline{\mathcal{L}}(\mathbb{Z}^{\prime}))={}_{-1}\overline{L_{0}%
}(\mathbb{Z}^{\prime})$ is extracted from the homotopy cartesian diagram%
\[%
\begin{array}
[c]{ccc}%
_{-1}\overline{\mathcal{L}}(\mathbb{Z}^{\prime}) & \longrightarrow &
_{-1}\mathcal{L}(\mathbb{R})\\
\downarrow &  & \downarrow\\
_{-1}\mathcal{L}(\mathbb{F}_{3}) & \longrightarrow & _{-1}\mathcal{L}%
(\mathbb{C})
\end{array}
\]
and reduces to the following:
\[
_{-1}L_{1}^{\mathrm{top}}(\mathbb{C)}\longrightarrow{}_{-1}\overline{L_{0}%
}(\mathbb{Z}^{\prime})\longrightarrow{}_{-1}L_{0}(\mathbb{F}_{3})\oplus{}%
_{-1}L_{0}(\mathbb{R})\twoheadrightarrow{}_{-1}L_{0}(\mathbb{C})\text{.}%
\]
Since
\[
_{-1}L_{1}^{\mathrm{top}}(\mathbb{C})=\pi_{1}(B\mathrm{Sp})=\pi_{5}%
(BO)=0\text{,}%
\]
this sequence shows that ${}_{-1}\overline{L_{0}}(\mathbb{Z}^{\prime})$ is the
same pull-back as $_{-1}L_{0}(\mathbb{Z}^{\prime})$, making $\pi_{0}({}%
_{-1}\phi)$ an isomorphism.

\begin{Void}
\label{-1L1}\textbf{The map} ${}_{-1}\phi:{}_{-1}\mathcal{L}(\mathbb{Z}%
^{\prime})_{\#}\longrightarrow{}_{-1}\overline{\mathcal{L}}(\mathbb{Z}%
^{\prime})_{\#}$ \textbf{induces an isomorphism on} $\pi_{1}$.
\end{Void}

After \cite{K:LNM343}p.382, we have
\[
_{-1}L_{1}(\mathbb{Z}^{\prime})=\mathrm{Sp}(\mathbb{Z}^{\prime})\left/
[\mathrm{Sp}(\mathbb{Z}^{\prime}),\,\mathrm{Sp}(\mathbb{Z}^{\prime})]\right.
=0\text{.}%
\]
On the other hand, the same homotopy cartesian square as above gives rise to
the following exact sequence (after completion)${}$%
\[
_{-1}L_{2}^{\mathrm{top}}(\mathbb{C})\longrightarrow{}_{-1}\overline{L}%
_{1}(\mathbb{Z}^{\prime})_{\#}\longrightarrow{}_{-1}L_{1}(\mathbb{F}_{3}%
)_{\#}\,\oplus\;{}_{-1}L_{1}^{\mathrm{top}}(\mathbb{R})_{\#}\longrightarrow
{}_{-1}L_{1}^{\mathrm{top}}(\mathbb{C})_{\#}%
\]
We have here, from ${}_{-1}O(\mathbb{C})\simeq\mathrm{Sp}$ and $_{-1}%
O(\mathbb{R})\simeq U$,%
\[
_{-1}L_{2}^{\mathrm{top}}(\mathbb{C})=\pi_{2}(B\mathrm{Sp})=\pi_{6}%
(BO)=0\text{;}%
\]%
\[
_{-1}L_{1}^{\mathrm{top}}(\mathbb{R})=\pi_{1}(BU)=0\text{;}%
\]%
\[
_{-1}L_{1}(\mathbb{F}_{3})=\mathrm{Sp}(\mathbb{F}_{3})\left/  \left[
\mathrm{Sp}(\mathbb{F}_{3}),\,\mathrm{Sp}(\mathbb{F}_{3})\right]  \right.
=0\text{.}%
\]
Therefore, $_{-1}\overline{L_{1}}(\mathbb{Z}^{\prime})$ is also reduced to $0$.

\begin{Void}
\label{1L0}\textbf{The map} ${}_{1}\phi:{}_{1}\mathcal{L}(\mathbb{Z}^{\prime
})_{\#}\longrightarrow{}_{1}\overline{\mathcal{L}}(\mathbb{Z}^{\prime})_{\#}$
\textbf{induces an isomorphism on} $\pi_{0}$.
\end{Void}

It is well-known that the Witt-Grothendieck group of $\mathbb{Z}^{\prime}$ is
isomorphic to $\mathbb{Z\oplus Z\oplus Z}/2$ (signature, rank and
discriminant; note that the discriminant is positive on real hyperbolic forms).

For the determination of ${}_{1}\overline{L_{0}}(\mathbb{Z}^{\prime})$, we
observe that the homotopy cartesian diagram%
\[%
\begin{array}
[c]{ccc}%
_{1}\overline{\mathcal{L}}(\mathbb{Z}^{\prime}) & \longrightarrow &
_{1}\mathcal{L}(\mathbb{R})\\
\downarrow &  & \downarrow\\
_{1}\mathcal{L}(\mathbb{F}_{3}) & \longrightarrow & _{1}\mathcal{L}%
(\mathbb{C})
\end{array}
\]
gives rise to the following exact sequence:{\small
\[
_{1}L_{1}(\mathbb{F}_{3})\oplus{}_{1}L_{1}^{\mathrm{top}}(\mathbb{R}%
)\rightarrow{}_{1}L_{1}^{\mathrm{top}}(\mathbb{C})\rightarrow{}_{1}%
\overline{L_{0}}(\mathbb{Z}^{\prime})\rightarrow{}_{1}L_{0}(\mathbb{F}%
_{3})\oplus\,{}_{1}L_{0}(\mathbb{R})\twoheadrightarrow{}_{1}L_{0}(\mathbb{C})
\]
}Recalling that $_{1}O(\mathbb{C})\simeq O$ and $_{1}O(\mathbb{R})\simeq
O\times O$, we have
\[
_{1}L_{1}^{\mathrm{top}}(\mathbb{C})\cong K_{1}^{\mathrm{top}}(\mathbb{R)}%
\cong\mathbb{Z}/2
\]
and
\[
_{1}L_{1}^{\mathrm{top}}(\mathbb{R})\cong K_{1}^{\mathrm{top}}(\mathbb{R)}%
\oplus K_{1}^{\mathrm{top}}(\mathbb{R})\cong\mathbb{Z}/2\oplus\mathbb{Z}%
/2\text{,}%
\]
and the first arrow of the above sequence is surjective. On the other hand, we
have the well-known isomorphisms%
\[%
\begin{array}
[c]{lll}%
_{1}L_{0}(\mathbb{C})\cong\mathbb{Z} & \qquad & \text{\cite{K:InvM73}p.248}\\
_{1}L_{0}(\mathbb{R})\cong\mathbb{Z\oplus Z} &  & \text{\cite{K:LNM343}%
p.306}\\
_{1}L_{0}(\mathbb{F}_{3})\cong\mathbb{Z\oplus Z}/2 &  & \text{\cite{Serre}%
IV.1.7}%
\end{array}
\]
which show that $_{1}\overline{L}_{0}(\mathbb{Z}^{\prime})$ is also isomorphic
to $\mathbb{Z\oplus Z}$ $\oplus\mathbb{\ Z}/2$.

It remains to check that the given map between these two isomorphic groups is
itself an isomorphism. The generator of the summand $\mathbb{Z}/2$ in
$_{1}L_{0}(\mathbb{Z}^{\prime})$ is the difference of the two quadratic forms
$\left\langle 1\right\rangle -\left\langle 2\right\rangle $ with the usual
notations. This element maps into the generator of the summand $\mathbb{Z}/2$
in $_{1}L_{0}(\mathbb{F}_{3})$ which is detected by the discriminant in
$(\mathbb{F}_{3}^{\bullet})/(\mathbb{F}_{3}^{\bullet})^{2}$. This discriminant
is just the Legendre symbol (for $p=3$)%
\[%
\genfrac{(}{)}{}{}{2}{p}%
=(-1)^{(p^{2}-1)/8}%
\]
which is equal to $-1$ (more simply, $2$ is not a square mod $3$). The
torsion-free summands correspond to the rank and index/signature invariants
\cite{Serre}IV, V, which are preserved respectively by the ring maps and the
inclusion $\mathbb{Z}^{\prime}\hookrightarrow\mathbb{R}$. Thus we obtain the
required isomorphism from $_{1}L_{0}(\mathbb{Z}^{\prime})\cong\mathbb{Z\oplus
Z\oplus Z}/2$ to the kernel ${}_{1}\overline{L}_{0}(\mathbb{Z}^{\prime})$ of
${}_{1}L_{0}(\mathbb{F}_{3})\oplus{}_{1}L_{0}(\mathbb{R})\twoheadrightarrow
{}_{1}L_{0}(\mathbb{C})$.

\begin{Void}
\label{1L1}\textbf{The map} ${}_{1}\phi:{}_{1}\mathcal{L}(\mathbb{Z}^{\prime
})_{\#}\longrightarrow{}_{1}\overline{\mathcal{L}}(\mathbb{Z}^{\prime})_{\#}$
\textbf{induces an isomorphism on} $\pi_{1}$.%
\nopagebreak
\end{Void}

From \cite{Bass: AmJM74}(4.7.6), there is an exact sequence%
\[
\mathrm{SK}_{1}(\mathbb{Z}^{\prime})\overset{H}{\longrightarrow}{}_{1}%
L_{1}(\mathbb{Z}^{\prime})\overset{(\mathrm{\det},\,\mathrm{SN})}%
{\longrightarrow}\mathrm{Ip}(\mathbb{Z}^{\prime})\oplus\mathrm{Discr}%
(\mathbb{Z}^{\prime})\rightarrow0\text{,}%
\]
where $\mathrm{Ip}(\mathbb{Z}^{\prime})=\{0_{\mathrm{Ip}\mathbb{Z}^{\prime}%
},\,1_{\mathrm{Ip}\mathbb{Z}^{\prime}}\}$ denotes the group of idempotents,
and there is a short exact sequence \cite{Bass: AmJM74}p.156%
\[
\mathbb{Z}^{\prime\bullet}/(\mathbb{Z}^{\prime\bullet})^{2}\rightarrowtail
\mathrm{Discr}(\mathbb{Z}^{\prime})\twoheadrightarrow{}_{2}\mathrm{Pic}%
(\mathbb{Z}^{\prime})\text{.}%
\]
Now $\mathrm{Pic}(\mathbb{Z}^{\prime})=0$ since $\mathbb{Z}^{\prime}$ is a
principal ideal domain. Further, from the localization sequence%
\[
K_{1}(\mathbb{F}_{2})\rightarrow K_{1}(\mathbb{Z})\rightarrow K_{1}%
(\mathbb{Z}^{\prime})\rightarrow K_{0}(\mathbb{F}_{2})\rightarrow
K_{0}(\mathbb{Z})\rightarrow K_{0}(\mathbb{Z}^{\prime})\text{,}%
\]
namely
\[
0\rightarrow\mathbb{Z}/2\rightarrow K_{1}(\mathbb{Z}^{\prime})\rightarrow
\mathbb{Z}\rightarrow\mathbb{Z}\rightarrow\mathbb{Z}\text{,}%
\]
we have that $K_{1}(\mathbb{Z}^{\prime})\cong\mathbb{Z}/2\oplus\mathbb{Z}%
\cong\mathbb{Z}^{\prime}{}^{\bullet}$. So $\mathrm{SK}_{1}(\mathbb{Z}^{\prime
})=0$. Hence we have%
\[
_{1}L_{1}(\mathbb{Z}^{\prime})\cong\mathbb{Z}/2\oplus(\mathbb{Z}^{\prime
})^{\bullet}/(\mathbb{Z}^{\prime\bullet})^{2}\cong\mathbb{Z}/2\oplus
\mathbb{Z}/2\oplus\mathbb{Z}/2
\]
generated respectively by $1_{\mathrm{Ip}\mathbb{Z}^{\prime}}$, $2(\mathbb{Z}%
^{\prime\bullet})^{2}$ and $(-1)(\mathbb{Z}^{\prime\bullet})^{2}$.

On the other hand, from the homotopy cartesian diagram above we deduce the
exact sequence{\small
\[
_{1}L_{2}(\mathbb{F}_{3})\oplus\,{}_{1}L_{2}^{\mathrm{top}}%
(\mathbb{R)\rightarrow{}}_{1}L_{2}^{\mathrm{top}}(\mathbb{C})\rightarrow{}%
_{1}\overline{L}_{1}(\mathbb{Z}^{\prime})\rightarrow\mathbb{{}}_{1}%
L_{1}(\mathbb{F}_{3})\oplus\,\mathbb{{}}_{1}L_{1}^{\mathrm{top}}%
(\mathbb{R})\,\mathbb{\rightarrow}_{1}L_{1}^{\mathrm{top}}(\mathbb{C})
\]
}Again, from the homotopy equivalences $_{1}O(\mathbb{R})\simeq O\times O$ and
$_{1}O(\mathbb{C})\simeq O$ it follows that the map
\[
\mathbb{Z}/2\oplus\mathbb{Z}/2\cong{}_{1}L_{2}^{\mathrm{top}}(\mathbb{R}%
)\longrightarrow{}_{1}L_{2}^{\mathrm{top}}(\mathbb{C})\cong\mathbb{Z}/2
\]
is surjective. Therefore, ${}_{1}\overline{L}_{1}(\mathbb{Z}^{\prime})$ is
identified with the kernel of the map%
\[
\mathbb{Z}/2\oplus\mathbb{Z}/2\oplus\mathbb{Z}/2\oplus\mathbb{Z}/2\cong{}%
_{1}L_{1}(\mathbb{F}_{3})\oplus{}_{1}L_{1}^{\mathrm{top}}(\mathbb{R}%
)\longrightarrow{}_{1}L_{1}^{\mathrm{top}}(\mathbb{C)}\cong\mathbb{Z}/2
\]
which is also the sum of three copies of $\mathbb{Z}/2$.

It remains to check that {}$_{1}\phi$ yields this isomorphism. First recall
from \cite{KV2}(5.6) that for a discrete field $F$ we have
\[
{}_{1}L_{1}(F)\cong\mathrm{Ip}(F)\oplus\mathrm{Discr}(F)\cong\mathbb{Z}%
/2\oplus F^{\bullet}/(F^{\bullet})^{2}\text{.}%
\]
In particular, this gives canonical generators for the cases of real and
complex fields, where the algebraic and topological ${}_{1}L_{1}$ groups
coincide. Then, by naturality of $\mathrm{Ip}$ and $\mathrm{Discr}$, we have
the three given generators of ${}_{1}L_{1}(\mathbb{Z}^{\prime})$ mapping
respectively under {}$_{1}\phi$ to the generators%
\[
1_{\mathrm{Ip}\mathbb{F}_{3}}+1_{\mathrm{Ip}\mathbb{R}}\text{,\quad
}2(\mathbb{F}_{3}^{\bullet})^{2}+(-1)(\mathbb{R}^{\bullet})^{2}\text{,\quad
}(-1)(\mathbb{R}^{\bullet})^{2}%
\]
of
\[
\mathrm{Ker}\left[  \left\langle 1_{\mathrm{Ip}\mathbb{F}_{3}}\right\rangle
\oplus\left\langle 2(\mathbb{F}_{3}^{\bullet})^{2}\right\rangle \oplus
\left\langle 1_{\mathrm{Ip}\mathbb{R}}\right\rangle \oplus\left\langle
(-1)(\mathbb{R}^{\bullet})^{2}\right\rangle \longrightarrow\left\langle
1_{\mathrm{Ip}\mathbb{C}}\right\rangle \right]  \text{.}%
\]
So we have the desired isomorphism after all.

\begin{Void}
\textbf{Summary: the} mod $2$ \textbf{equivalence for} $\pi_{0}$, $\pi_{1}$%
\nopagebreak
\end{Void}

We now use the data collected above to display the low-dimensional information
in the form we need.

\begin{Lemma}
For $\varepsilon\in\{\pm1\},$ $i\in\{0,1\}$, the map ${}_{\varepsilon}\phi
:{}_{\varepsilon}\mathcal{L}(\mathbb{Z}^{\prime})_{\#}\longrightarrow
{}_{\varepsilon}\overline{\mathcal{L}}(\mathbb{Z}^{\prime})_{\#}$ induces an
isomorphism on $\pi_{i}$ with $\mathbb{Z}/2$ coefficients.
\end{Lemma}

\noindent\textbf{Proof.\quad} Since we already have this result with integer
coefficients, to obtain the mod $2$ version via the universal coefficient
sequence, it suffices to check that there is induced an integral isomorphism
on $\pi_{-1}$ of these spectra. Here the spaces $B_{\varepsilon}O(\mathbb{R}%
)$, $B_{\varepsilon}O(\mathbb{C})$ are expressible in terms of $BO$, $BU$ and
$B\mathrm{Sp}$, as noted above, and so contribute zero groups. Therefore we
need only check that ${}_{\varepsilon}\phi$ induces isomorphisms from
${}_{\varepsilon}L_{-1}(\mathbb{Z}^{\prime})$ to ${}_{\varepsilon}%
L_{-1}(\mathbb{F}_{3})$. To do this, essentially in a reprise of downward
induction in Theorem \ref{K-L bootstrap}, we exploit the natural exact
sequences (for $A=\mathbb{Z}^{\prime},\mathbb{F}_{3}$)%
\[
K_{1}(A)\rightarrow{}_{-\varepsilon}L_{1}(A)\rightarrow{}_{-\varepsilon}%
U_{0}(A)\rightarrow K_{0}(A)\rightarrow{}_{-\varepsilon}L_{0}(A)
\]
and (via the fundamental theorem \cite{K:AnnM112 fun})%
\[
K_{0}(A)\rightarrow{}_{-\varepsilon}U_{0}(A)\rightarrow{}_{\varepsilon}%
L_{-1}(A)\rightarrow K_{-1}(A)
\]
We of course use the information we've already gathered about these groups.

\smallskip

\noindent$\varepsilon=1$\textbf{.\quad} For both $A=\mathbb{Z}^{\prime
},\mathbb{F}_{3}$ the first sequence terminates with
\[
0\rightarrow{}_{-1}U_{0}(A)\rightarrow\mathbb{Z\rightarrow Z}%
\]
where the final map is given by the rank (see (\ref{-1L0})) and so is nonzero.
So ${}_{-1}U_{0}(A)$ vanishes, and the second sequence finishes with
\[
0\rightarrow{}_{1}L_{-1}(A)\rightarrow0\text{,}%
\]
leaving ${}_{1}L_{-1}(\mathbb{Z}^{\prime})={}_{1}L_{-1}(\mathbb{F}_{3})=0$.

\smallskip

\noindent$\varepsilon=-1$\textbf{.\quad} Here the first map of exact sequences
becomes%
\[%
\begin{array}
[c]{ccccccccc}%
\mathbb{Z\oplus Z}/2 & \rightarrow & (\mathbb{Z}/2)^{3} & \rightarrow & {}%
_{1}U_{0}(\mathbb{Z}^{\prime}) & \rightarrow & \mathbb{Z} & \rightarrow &
\mathbb{Z}^{2}\mathbb{\oplus Z}/2\\
\downarrow &  & \downarrow &  & \downarrow &  & \downarrow &  & \downarrow\\
\mathbb{Z}/2 & \rightarrow & (\mathbb{Z}/2)^{2} & \rightarrow & {}_{1}%
U_{0}(\mathbb{F}_{3}) & \rightarrow & \mathbb{Z} & \rightarrow &
\mathbb{Z\oplus Z}/2
\end{array}
\]
where in each sequence the final map is injective and the first has cokernel
$\mathbb{Z}/2$, by our previous considerations. We therefore obtain that
${}_{\varepsilon}\phi$ induces an isomorphism ${}_{1}U_{0}(\mathbb{Z}^{\prime
})\overset{\cong}{\longrightarrow}{}_{1}U_{0}(\mathbb{F}_{3})\cong
\mathbb{Z}/2$. Since $K_{i}(\mathbb{Z}^{\prime})\rightarrow K_{i}%
(\mathbb{F}_{3})$ is an isomorphism for $i=0,-1$, the desired isomorphism on
${}_{-1}L_{-1}(A)$ now follows by applying the five lemma to the second
sequence.\hfill$\Box\medskip$

\section{\label{htpy equiv Thm A}The homotopy equivalence for Theorem A}

In accordance with the argument of \cite{Sullivan}p.32 (exploiting Bockstein
homomorphisms), our claimed homotopy equivalence%
\[
{}_{\varepsilon}\phi:{}_{\varepsilon}\mathcal{L}(\mathbb{Z}^{\prime}%
)_{\#}\longrightarrow{}_{\varepsilon}\overline{\mathcal{L}}(\mathbb{Z}%
^{\prime})_{\#}%
\]
of $2$-completed spaces may be deduced by proving that the map is both a
rational and a $2$-local equivalence. Thus our argument breaks into two parts.

\medskip\newpage

\begin{Void}
\textbf{We work rationally}.%
\nopagebreak
\end{Void}

Since by \cite{F} the space $_{\varepsilon}\mathcal{L}(\mathbb{F}_{3})_{\#}$
is rationally $\mathbb{Q}_{\#}$, the space ${}_{\varepsilon}\overline
{\mathcal{L}}(\mathbb{Z}^{\prime})_{\#}$ rationally reduces to the homotopy
fiber of the rationalized map from ${}_{\varepsilon}\mathcal{L}(\mathbb{R}%
)_{\#}$ to ${}_{\varepsilon}\mathcal{L}(\mathbb{C})_{\#}$, namely
$\mathbb{Q}_{\#}\times\Omega^{7+\varepsilon}(BO)$. Now by \cite{Bor}\S12, the
map from ${}_{\varepsilon}\mathcal{L}(\mathbb{Z}^{\prime})_{\#}$ to
$\mathbb{Q}_{\#}\times\Omega^{7+\varepsilon}(BO)$ is a rational homotopy
equivalence. (Recall that finite generation of the homotopy groups of
${}_{\varepsilon}\mathcal{L}(\mathbb{Z}^{\prime})$ was established in
Corollary \ref{fg}.) For the situation with $\mathbb{Z}^{\prime}$
coefficients, see \cite{K:AnnM112 Quillen}p.253 et seq.

\begin{Void}
\textbf{We work} $\mathrm{mod\,}2$ (\textbf{which is prime to} $3$).
\end{Void}

In fact, for the arguments we present, it is no more difficult to work in the
more general setting mod\ $\ell$, and consider the finite field $\mathbb{F}%
_{q}$ with $q$ elements, $q$ odd, and $\ell$ coprime to $q$. We require some
analysis of Quillen's Brauer lifting $\mathcal{K}(\mathbb{F}_{q}%
)\longrightarrow$ $\mathcal{K}^{\mathrm{top}}(\mathbb{C})$ \cite{Q} and its
refinement by E.\thinspace M.~ Friedlander \cite{F} in the Hermitian case as a
map $_{\varepsilon}\mathcal{L}(\mathbb{F}_{q})\longrightarrow{}_{\varepsilon
}\mathcal{L}^{\mathrm{top}}(\mathbb{C})$ (see Theorem \ref{Friedlander fibn}
below). By the theorem of Quillen we have a homotopy fibration%
\begin{equation}%
\begin{array}
[c]{ccccc}%
\mathcal{K}(\mathbb{F}_{q}) & \longrightarrow & \mathcal{K}^{\mathrm{top}%
}(\mathbb{C}) & \longrightarrow & \mathcal{K}^{\mathrm{top}}(\mathbb{C)}%
\end{array}
\label{QF fibrations}%
\end{equation}
where the second arrow is defined by $\psi^{q}-1$, $\psi^{q}$ being the Adams
operation. In other words, we have a homotopy equivalence between
$\mathcal{K}(\mathbb{F}_{q})$ and a certain homotopy fiber.

Our strategy is to replace the map ${}_{\varepsilon}\mathcal{L}(\mathbb{F}%
_{q})\longrightarrow{}_{\varepsilon}\mathcal{L}^{\mathrm{top}}(\mathbb{C})$ by
one induced from a map of rings, so that we obtain a ring $B$ whose
${}_{\varepsilon}\mathcal{L}(B)$ serves as a candidate for ${}_{\varepsilon
}\overline{\mathcal{L}}(\mathbb{Z}^{\prime})$. To do this, we of course take
advantage of the fact that we are interested only in homotopy equivalences
$\operatorname{mod}$\ $\ell$. Thus, for example, by \cite{S} it does not
matter whether we use the usual or discrete topology on $\mathbb{R}$ and
$\mathbb{C}$. In order to fix ideas (the usual topology leads to alternative
arguments), we take the discrete topology. We proceed in five steps.

First, we exploit the fact \cite{S} that in $\operatorname{mod}$\ $\ell$
$K$-theory we may replace $\mathbb{F}_{q}$ by the ring $\mathbb{Z}_{q}$ of
$q$-adic integers, the Witt ring of $\mathbb{F}_{q}$. Likewise, $\mathbb{\bar
{Z}}_{q}$ is the ring of Witt vectors of the algebraic closure $\mathbb{\bar
{F}}_{q}$ of $\mathbb{F}_{q}$, with $\mathbb{\bar{Q}}_{q}$ is the field of
fractions of $\mathbb{\bar{Z}}_{q}$.

\begin{Lemma}
\label{Zq and Fq}The homomorphism $\mathbb{Z}_{q}\rightarrow\mathbb{F}_{q}$
induces an isomorphism of all $\operatorname{mod}$\ $\ell$ $L$-groups.
\end{Lemma}

\noindent\textbf{Proof.\quad} According to Corollary \ref{L to mod p L} above,
we just have to check that $\mathbb{Z}_{q}\rightarrow\mathbb{F}_{q}$ has the
following effects:

\noindent\textbf{Isomorphism on }${}_{\varepsilon}L_{0}$. This follows from
\cite{MH}p.7 when $\varepsilon=1$ (both groups $\mathbb{Z}$), and from
\cite{BakScharlau} when $\varepsilon=-1$ (both groups $\mathbb{Z\oplus Z}/2$).

\noindent\textbf{Isomorphism on }${}_{\varepsilon}L_{1}$. For $\varepsilon=1$
($\mathbb{Z}/2\mathbb{\oplus Z}/2$), see \cite{Bass: AmJM74}, and for
$\varepsilon=-1$ (both groups vanish), see \cite{BassMilnorSerre}(13.2).

\noindent\textbf{Isomorphism on }$K_{-1}$. Both groups vanish since the rings
are regular Noetherian.

\noindent\textbf{Isomorphism on }$K_{0}$. Since each ring is local, the groups
are $\mathbb{Z}$.

\noindent\textbf{Epimorphism on }$K_{1}$. We have $K_{1}(\mathbb{Z}%
_{q})=\mathrm{GL}_{1}(\mathbb{Z}_{q})\twoheadrightarrow\mathrm{GL}%
_{1}(\mathbb{F}_{q})=K_{1}(\mathbb{F}_{q})$.\hfill$\Box\medskip$

Second, we observe that in $K$-theory with finite coefficients the Adams
operation $\psi^{q}$ on $\mathcal{K}(\mathbb{C})$ is induced by a certain
(non-unique) discontinuous automorphism of $\mathbb{C}$ described as follows.
First, we choose an embedding of the $q$-adic numbers $\mathbb{Q}_{q}$ in
$\mathbb{C}$, and consider the union of all the cyclotomic fields (over
$\mathbb{Q}_{q}$) generated by the $n\,$th roots of unity where $n=q^{r}-1$
with $r$ a power of $q$. This field admits an automorphism induced by the map
sending a primitive root $\zeta$ to $\zeta^{q}$. After Zorn's lemma, the
automorphism extends to a (discontinuous) automorphism $\psi$ of $\mathbb{C}$
that induces $\psi^{q}$ \cite{Quillen: Top68}. The homotopy fiber of $\psi
^{q}-1:\mathcal{K}(\mathbb{C})\rightarrow\mathcal{K}(\mathbb{C})$ may now be
realized as $\mathcal{K}(F)$ where the equalizer ring $F$ of $\psi$ and
$\mathrm{id}$ is constructed in Appendix B as the pull-back%
\[%
\begin{array}
[c]{ccc}%
F & \twoheadrightarrow & \mathbb{C}\\
\downarrow & {}^{\cdot}\lrcorner & \quad\downarrow^{(\psi,\mathrm{id})}\\
M(\Delta) & \overset{\tilde{\Delta}}{\twoheadrightarrow} & \mathbb{C}%
\times\mathbb{C}%
\end{array}
\]

Next, we use the fact that $\mathbb{Z}_{q}$ is fixed by $\psi^{q}$, in order
to factorize its inclusion in $\mathbb{C}$ as
\[
\mathbb{Z}_{q}\longrightarrow F\overset{\eta}{\twoheadrightarrow}%
\mathbb{C}\text{.}%
\]
The induced map $\mathcal{K}(\mathbb{Z}_{q})\rightarrow\mathcal{K}(F)$ is then
a mod\ $\ell$ homotopy equivalence, by the upper part of the folklore
commuting diagram ($q$ odd)%
\begin{equation}%
\begin{array}
[c]{ccccc}%
\mathcal{K}(\mathbb{F}_{q}) & \longrightarrow & \mathcal{K}^{\mathrm{top}%
}(\mathbb{C}) & \overset{\psi^{q}-1}{\longrightarrow} & \mathcal{K}%
^{\mathrm{top}}(\mathbb{C})\\
\uparrow &  & \uparrow &  & \uparrow\\
\mathcal{K}(\mathbb{Z}_{q}) & \longrightarrow & \mathcal{K}(\mathbb{\bar{Q}%
}_{q}) & \overset{\psi-1}{\longrightarrow} & \mathcal{K}(\mathbb{\bar{Q}}%
_{q})\\
&  & \uparrow &  & \uparrow\\
\downarrow &  & \mathcal{K}(\mathbb{\bar{Z}}_{q}) & \overset{\psi
-1}{\longrightarrow} & \mathcal{K}(\mathbb{\bar{Z}}_{q})\\
&  & \downarrow &  & \downarrow\\
\mathcal{K}(\mathbb{F}_{q}) & \longrightarrow & \mathcal{K}(\mathbb{\bar{F}%
}_{q}) & \overset{\psi-1}{\longrightarrow} & \mathcal{K}(\mathbb{\bar{F}}_{q})
\end{array}
\label{folklore}%
\end{equation}
of Quillen \cite{Quillen: AnnM cohom gps}, Friedlander and Suslin \cite{S}, in
which all vertical arrows, induced by ring homomorphisms, are mod\ $\ell$
equivalences ($\ell$ prime to $q$), and $\psi$ is the Frobenius map.

Fourth, we consider the commuting diagram of ring homomorphisms, in which the
two right-hand squares are pull-backs, the bottom composite is just the
diagonal map $\Delta:\mathbb{C}\rightarrow\mathbb{C}\times\mathbb{C}$, and the
maps $\varphi,\varphi^{\prime}$ derive from the fact that $\psi$ fixes
$\mathbb{Z}_{q}$.%
\[%
\begin{array}
[c]{ccccc}%
\mathbb{Z}^{\prime} & \overset{\varphi^{\prime}}{\longrightarrow} & F^{\prime}%
& \twoheadrightarrow & \mathbb{R}\\
\downarrow &  & \downarrow & {}^{\cdot}\lrcorner & \quad\downarrow
^{\mathrm{inc}}\\
\mathbb{Z}_{q} & \overset{\varphi}{\longrightarrow} & F & \overset{\eta
}{\twoheadrightarrow} & \mathbb{C}\\
\downarrow &  & \downarrow & {}^{\cdot}\lrcorner & \quad\downarrow
^{(\psi,\mathrm{id})}\\
\mathbb{C} & \longrightarrow &  M(\Delta) & \overset{\tilde{\Delta}%
}{\twoheadrightarrow} & \mathbb{C}\times\mathbb{C}%
\end{array}
\]

Finally, from Theorem \ref{M-V} above we have the following.

\begin{Lemma}
There are $\mathrm{mod\,}\ell$ homotopy equivalences between $\mathcal{K}%
(F^{\prime})$ and $\mathcal{\bar{K}(}\mathbb{Z}^{\prime})$, and between
${}_{\varepsilon}\mathcal{L}(F^{\prime})$ and ${}_{\varepsilon}\overline
{\mathcal{L}}(\mathbb{Z}^{\prime})$.\hfill$\Box$
\end{Lemma}

We remark that alternative proofs of the lemma above may be obtained from the
fact, as in \cite{Mitchell: plus}, that the inclusion $\mathbb{Z}%
_{q}\hookrightarrow\mathbb{C}$ induces Brauer lifting. One can proceed with
the discrete topology on $\mathbb{C}$, as above, or, with the usual topology
on $\mathbb{C}$ by replacing $\mathbb{C}[x]$ by the Banach ring $\mathbb{C}%
\left\langle x\right\rangle $ of complex convergent power series, as in
\cite{KV1}.

To finish the proof of Theorem A, it now remains to compare ${}_{\varepsilon
}\mathcal{L}(\mathbb{Z}^{\prime})$ and ${}_{\varepsilon}\mathcal{L}(F^{\prime
})$. By combining the above lemma with the equivalence $\mathcal{K}%
(\mathbb{Z}^{\prime})_{\#}\rightarrow\mathcal{\bar{K}}(\mathbb{Z}^{\prime
})_{\#}$ discussed at the beginning of the paper, we have that the ring
homomorphism $\varphi^{\prime}:\mathbb{Z}^{\prime}\rightarrow F^{\prime}$
induces a $\mathrm{mod\,}2$ homotopy equivalence $\mathcal{K}(\mathbb{Z}%
^{\prime})\rightarrow\mathcal{K}(F^{\prime})$. From Section \ref{low}, the
induced map ${}_{\varepsilon}\mathcal{L}(\mathbb{Z}^{\prime})\rightarrow
{}_{\varepsilon}\mathcal{L}(F^{\prime})$ also induces an isomorphism on the
homotopy groups $\pi_{0}$ and $\pi_{1}$ $\mathrm{mod\,}2$. Hence Theorem A is
now a consequence of Theorem \ref{K-L bootstrap}.\hfill$\Box\smallskip$

With the methods developed here, we can also prove the $L$-theory counterpart
of the fibration (\ref{QF fibrations}), a slightly more general formulation of
Friedlander's result \cite{F}, \cite{FiedorowiczPriddy}, to be used in Lemma
\ref{hZ/2 for F_3}.

\begin{Theorem}
\label{Friedlander fibn}For any odd prime power $q$, there are homotopy
fibrations of spaces
\[%
\begin{array}
[c]{ccccc}%
{}_{\varepsilon}\mathcal{L}(\mathbb{F}_{q}) & \longrightarrow & {}%
_{\varepsilon}\mathcal{L}^{\mathrm{top}}(\mathbb{C}) & \longrightarrow &
{}_{\varepsilon}\mathcal{L}^{\mathrm{top}}(\mathbb{C)}\text{,}%
\end{array}
\]
where the second arrow is defined by $\psi^{q}-1$, and%
\[%
\begin{array}
[c]{ccccc}%
{}_{\varepsilon}\mathcal{L}(\mathbb{F}_{q}) & \longrightarrow & {}%
_{\varepsilon}\mathcal{L}(\mathbb{\bar{F}}_{q}) & \longrightarrow &
{}_{\varepsilon}\mathcal{L}(\mathbb{\bar{F}}_{q}\mathbb{)}\text{,}%
\end{array}
\]
with the second arrow defined by $\psi-1$ (and $\psi$ the Frobenius automorphism).
\end{Theorem}

\noindent\textbf{Proof.\quad} We are essentially using Friedlander's
computations \cite{F} of ${}_{1}L(\mathbb{F}_{q})$ for $i=5,6$ and of ${}%
_{-1}L(\mathbb{F}_{q})$ for $i=1,2$.

First note that the theorem holds both rationally and modulo the
characteristic of $\mathbb{F}_{q}$, since the action of $\psi^{q}-1$ on
$\pi_{2n}$ is as $q^{n}-1$, and so invertible. Therefore we now work
mod\ $\ell$ where $\ell$ is prime to $q$.

For the first homotopy fibration, we can repeat the argument used above, but
now with respect to a map instead of $\varphi^{\prime}:\mathbb{Z}^{\prime
}\rightarrow F^{\prime}$. In other words, we write the homotopy fiber of
${}_{\varepsilon}\mathcal{L}^{\mathrm{top}}(\mathbb{C})\overset{\psi^{q}%
-1}{\longrightarrow}{}_{\varepsilon}\mathcal{L}^{\mathrm{top}}(\mathbb{C})$ as
${}_{\varepsilon}\mathcal{L}(F)$, and we have a ring map $\varphi
:\mathbb{Z}_{q}\rightarrow F$ that we want to show induces an isomorphism
${}_{\varepsilon}\tilde{L}_{i}(\mathbb{Z}_{q})\rightarrow{}_{\varepsilon
}\tilde{L}_{i}(F)$ of mod\ $\ell$ $L$-groups. By combining Lemma \ref{Zq and
Fq} and Friedlander's computation of ${}_{\varepsilon}\tilde{L}(\mathbb{F}%
_{q})$, we know that%
\begin{align*}
{}_{1}\tilde{L}_{i}(\mathbb{Z}_{q})  &  ={}_{1}\tilde{L}_{i}(F)=0\qquad
\text{for }i=5,6,\\
{}_{-1}\tilde{L}_{i}(\mathbb{Z}_{q})  &  ={}_{-1}\tilde{L}_{i}(F)=0\qquad
\text{for }i=1,2\text{.}%
\end{align*}
Therefore, using Theorem \ref{K-L bootstrap} in both directions, we have
\[
{}_{\varepsilon}\tilde{L}_{i}(\mathbb{Z}_{q})={}_{\varepsilon}\tilde{L}%
_{i}(F)\qquad\text{for }i\geq0\text{.}%
\]
(Note however that ${}_{\varepsilon}\tilde{L}_{-1}(F)\neq0={}_{\varepsilon
}\tilde{L}_{-1}(\mathbb{F}_{q})$.)

For the second fibration, we again use the folklore diagram (\ref{folklore}).
According to Appendix B, the homotopy fiber of the ${}_{\varepsilon
}\mathcal{L}$ space map induced by $\psi-1$ may be thought of as
${}_{\varepsilon}\mathcal{L}(B)$, where $B$ is the pull-back%
\[%
\begin{array}
[c]{ccc}%
B & \twoheadrightarrow & \mathbb{\bar{F}}_{q}\\
\downarrow & {}^{\cdot}\lrcorner & \quad\downarrow^{(\psi,\mathrm{id})}\\
M(\Delta) & \overset{\tilde{\Delta}}{\twoheadrightarrow} & \mathbb{\bar{F}%
}_{q}\times\mathbb{\bar{F}}_{q}%
\end{array}
\]
Because $\psi$ fixes $\mathbb{F}_{q}$, there is a canonical map $\mathbb{F}%
_{q}\longrightarrow B$ which we seek to show induces a mod\ $\ell$ homotopy
equivalence on ${}_{\varepsilon}\mathcal{L}$ spaces. After Quillen
\cite{Quillen: AnnM cohom gps}, this is known for the $\mathcal{K}$ spaces. So
the result follows from Corollary \ref{L to mod p L} as before, since the
$L$-theories of $\mathbb{\bar{F}}_{q}$ and $\mathbb{C}$ with finite
coefficients (coprime to $q$) coincide by \cite{K:InvM73}.\hfill
$\Box\smallskip$

\section{\label{L groups}Computations modulo odd torsion}

\begin{Theorem}
\label{Thm B in 6}Modulo a finite group of odd order, the groups
${}_{\varepsilon}L_{i}(\mathbb{Z}^{\prime})$ for $i\geq0$ are as follows,
where $\delta_{i0}$ denotes the Kronecker delta, and $2^{t}$ is the
$2$-primary part of $i+1$.
\[%
\begin{array}
[c]{ccc}%
i\;(\mathrm{\operatorname{mod}\,}8) & _{-1}L_{i}(\mathbb{Z}^{\prime}) &
_{1}L_{i}(\mathbb{Z}^{\prime})\\
0 & \delta_{i0}\mathbb{Z} & \delta_{i0}\mathbb{Z\oplus Z\oplus Z}/2\\
1 & 0 & \mathbb{Z}/2\oplus\mathbb{Z}/2\oplus\mathbb{Z}/2\\
2 & \mathbb{Z} & \mathbb{Z}/2\oplus\mathbb{Z}/2\\
3 & \mathbb{Z}/16 & \mathbb{Z}/8\\
4 & \mathbb{Z}/2 & \mathbb{Z}\\
5 & \mathbb{Z}/2 & 0\\
6 & \mathbb{Z} & 0\\
7 & \mathbb{Z}/2^{t+1} & \mathbb{Z}/2^{t+1}%
\end{array}
\]
Moreover, the odd torsion subgroup of $_{\varepsilon}L_{i}(\mathbb{Z}^{\prime
})$ is the invariant part of the odd torsion of $K_{i}(\mathbb{Z})$ induced by
the involution $M\mapsto{}^{\mathrm{t}}M^{-1}$ of $\mathrm{GL}(\mathbb{Z})$.
\end{Theorem}

\noindent\textbf{Proof.\quad} We first consider odd torsion. From the cases
$i=0,1$ dealt with above, neither $K_{i}(\mathbb{Z}^{\prime})$ nor
$_{\varepsilon}L_{i}(\mathbb{Z}^{\prime})$ contains (nontrivial) odd torsion,
so that the Witt and co-Witt groups
\begin{align*}
_{\varepsilon}W_{i}(\mathbb{Z}^{\prime})  &  =\mathrm{Coker\,}[K_{i}%
(\mathbb{Z}^{\prime})\longrightarrow{}_{\varepsilon}L_{i}(\mathbb{Z}^{\prime
})],\\
_{\varepsilon}W_{i}^{\prime}(\mathbb{Z}^{\prime})  &  =\mathrm{Ker\,}%
[{}_{\varepsilon}L_{i}(\mathbb{Z}^{\prime})\longrightarrow K_{i}%
(\mathbb{Z}^{\prime})]
\end{align*}
are odd-torsion-free for $\varepsilon\in\{\pm1\},$ $i\in\{0,1\}$. Then by
Proposition 1.35 and Theor\`{e}me 3.7 of \cite{K:AnnM112 Quillen} these groups
are odd-torsion-free for all $i$. In particular, the maps $K_{i}%
(\mathbb{Z}^{\prime})\longrightarrow{}_{\varepsilon}L_{i}(\mathbb{Z}^{\prime
})$ are surjective on odd torsion. The assertion about the odd torsion
subgroups then follows from the argument of \cite{K:AnnM112 Quillen}p.253
(specifically, using the weak equivalence there denoted $\mathcal{C}(A)\sim
{}_{\varepsilon}\mathcal{L}(A)$).

We can now work modulo odd torsion. The groups for $i=0$ have been discussed
above, so we may assume that $i\geq1$.

\smallskip

\noindent$\varepsilon=1$\textbf{.\quad} Since the homotopy cartesian diagram
\[%
\begin{array}
[c]{ccc}%
B_{1}O(\mathbb{Z}^{\prime})_{\#}^{+} & \longrightarrow &  B_{1}O(\mathbb{R}%
)_{\#}\\
\downarrow &  & \downarrow\\
B_{1}O(\mathbb{F}_{3})_{\#}^{+} & \longrightarrow &  B_{1}O(\mathbb{C})_{\#}%
\end{array}
\]
may also be written as
\[%
\begin{array}
[c]{ccc}%
B_{1}O(\mathbb{Z}^{\prime})_{\#}^{+} & \longrightarrow &  BO_{\#}\times
BO_{\#}\\
\downarrow &  & \downarrow\\
B_{1}O(\mathbb{F}_{3})_{\#}^{+} & \longrightarrow &  BO_{\#}%
\end{array}
\]
we have a split short exact sequence:
\[
0\rightarrow{}_{1}L_{i}(\mathbb{Z}^{\prime})\rightarrow\pi_{i}(BO)\oplus
\pi_{i}(BO)\oplus\pi_{i}(BO(\mathbb{F}_{3})^{+})\rightarrow\pi_{i}%
(BO)\rightarrow0
\]
which can simply be written as an isomorphism
\[
{}_{1}L_{i}(\mathbb{Z}^{\prime})\cong\pi_{i}(BO)\oplus\pi_{i}(BO(\mathbb{F}%
_{3})^{+})
\]
for $i>0$. (In fact, more precisely one has a homotopy decomposition
$B\,_{1}O(\mathbb{Z}^{\prime})^{+}\simeq BO\times B\,_{1}O(\mathbb{F}_{3}%
)^{+}$.) Then the results tabulated above follow immediately from the
computations of Friedlander \cite{F}, given the following well-known lemma,
easily proven by induction (or factorization). The relevance of $3^{r}-1$ here
is that it is the order of $K_{2r-1}(\mathbb{F}_{3})$.

\noindent\textbf{Numerical Claim.\quad} \emph{Write} $(r)_{2}$ \emph{for the}
$2$\emph{-primary part of} $r$\emph{. Then, for even }$r$\emph{,}
\[
\left(  3^{r}-1\right)  _{2}=4(r)_{2}\text{.}%
\]

\smallskip

\noindent$\varepsilon=-1$\textbf{.\quad} Using the homotopy cartesian square%
\[%
\begin{array}
[c]{ccc}%
B_{-1}O(\mathbb{Z}^{\prime})_{\#}^{+} & \longrightarrow &  B_{-1}%
O(\mathbb{R})_{\#}\\
\downarrow &  & \downarrow\\
B_{-1}O(\mathbb{F}_{3})_{\#}^{+} & \longrightarrow &  B_{-1}O(\mathbb{C})_{\#}%
\end{array}
\]
we argue as follows, with all spaces and homotopy groups $2$-completed. The
above reduces to the homotopy cartesian square%
\[%
\begin{array}
[c]{ccc}%
B\mathrm{Sp}\mathbb{Z}^{\prime}{}^{+} & \longrightarrow &  BU\\
\downarrow &  & \downarrow^{H}\\
B\mathrm{Sp}\mathbb{F}_{3}{}^{+} & \longrightarrow &  B\mathrm{Sp}%
\end{array}
\]
with vertical homotopy fiber $\mathrm{Sp}/U\simeq\Omega^{6}BO$ induced by the
hyperbolic map $H$, and horizontal homotopy fiber $\mathrm{Sp}\simeq\Omega
^{5}BO$. This gives rise to three exact homotopy sequences for each $i\geq
0:${\small
\[
\mathrm{A}_{i}:\qquad\pi_{i+1}B\mathrm{Sp}\mathbb{F}_{3}{}^{+}\rightarrow
\pi_{i+6}BO\rightarrow\pi_{i}B\mathrm{Sp}\mathbb{Z}^{\prime}{}^{+}%
\rightarrow\pi_{i}B\mathrm{Sp}\mathbb{F}_{3}{}^{+}\rightarrow\pi_{i+5}BO
\]%
\begin{align*}
\mathrm{B}_{i}  &  :\qquad\pi_{i+1}B\mathrm{Sp}\mathbb{F}_{3}{}^{+}\oplus
\pi_{i+1}BU\rightarrow\pi_{i+5}BO\rightarrow\pi_{i}B\mathrm{Sp}\mathbb{Z}%
^{\prime}{}^{+}\rightarrow\\
&  \rightarrow\pi_{i}B\mathrm{Sp}\mathbb{F}_{3}{}^{+}\oplus\pi_{i}%
BU\rightarrow\pi_{i+4}BO
\end{align*}%
\[
\mathrm{C}_{i}:\qquad\pi_{i+1}BU\rightarrow\pi_{i+5}BO\rightarrow\pi
_{i}B\mathrm{Sp}\mathbb{Z}^{\prime}{}^{+}\rightarrow\pi_{i}BU\rightarrow
\pi_{i+4}BO
\]
}

The groups $\pi_{i}B\mathrm{Sp}\mathbb{F}_{3}{}^{+}$, together with
information about the map $H_{\ast}:\pi_{i}BU\rightarrow\pi_{i+4}BO$ appearing
in $\mathrm{B}_{i}$, are calculated in \cite{F}. It follows immediately from
$\mathrm{A}_{i}$ that for $i\equiv0,1,2\;(\mathrm{\operatorname{mod}\,}8)$ the
map $\pi_{i+6}BO\rightarrow\pi_{i}B\mathrm{Sp}\mathbb{Z}^{\prime}{}^{+}$ is an
isomorphism, while for $i\equiv7\;(\mathrm{\operatorname{mod}\,}8)$ it is
$\pi_{i}B\mathrm{Sp}\mathbb{Z}^{\prime}{}^{+}\rightarrow\pi_{i}B\mathrm{Sp}%
\mathbb{F}_{3}{}^{+}$ that is an isomorphism (in this case one again appeals
to the above numerical claim). Likewise, when $i\equiv
6\;(\mathrm{\operatorname{mod}\,}8)$ the sequence $\mathrm{C}_{i}$ immediately
becomes:%
\[
0\rightarrow0\rightarrow\mathbb{Z}\rightarrow\mathbb{Z}\rightarrow
\mathbb{Z}/2
\]

For $i\equiv3\;(\mathrm{\operatorname{mod}\,}8)$, since $H_{\ast}$ is
multiplication by $2$ on $\pi_{i+1}BU\cong\pi_{i+5}BO\cong\mathbb{Z}$, from
$\mathrm{B}_{i}$ we have (using the numerical claim above) that $\pi
_{i}B\mathrm{Sp}\mathbb{Z}^{\prime}{}^{+}$ has order $16$. On the other hand,
from $\mathrm{C}_{i}$ this must be a cyclic group.

Then for $i\equiv4\;(\mathrm{\operatorname{mod}\,}8)$ the fact that $H_{\ast}$
is injective in $\mathrm{B}_{i}$ forces $\pi_{i}B\mathrm{Sp}\mathbb{Z}%
^{\prime}{}^{+}$ to be finite. Hence $\mathrm{C}_{i}$ reduces to the sequence%
\[
0\rightarrow\mathbb{Z}/2\rightarrow\mathbb{Z}/2\rightarrow\mathbb{Z}%
\rightarrowtail\mathbb{Z}%
\]
Since the sequence $\mathrm{A}_{i-1}$ began with the zero map, so
$\mathrm{A}_{i}$ ends with it, and becomes:%
\[
\mathbb{Z}/2\oplus\mathbb{Z}/2\twoheadrightarrow\mathbb{Z}/2\rightarrow
\mathbb{Z}/2\rightarrow\mathbb{Z}/2\overset{0}{\longrightarrow}\mathbb{Z}/2
\]

Finally, it follows that when $i\equiv5\;(\mathrm{\operatorname{mod}\,}8)$ the
final map of the sequence $\mathrm{A}_{i}$ (being the initial map of the
sequence $\mathrm{A}_{i-1}$) is an epimorphism with kernel $\pi_{i}%
B\mathrm{Sp}\mathbb{Z}^{\prime}{}^{+}=\mathbb{Z}/2$.\hfill$\Box$

\begin{Remarks}
\textbf{1. }Of course the $2$-primary torsion of $Z/2^{t+1}$ above is just the
$2$-primary torsion of Im $J$. Observe too that ${}_{1}L_{3}(\mathbb{Z}%
^{\prime})=\mathbb{Z}/24$, a less exotic result than $K_{3}(\mathbb{Z}%
^{\prime})$ which is $\mathbb{Z}/48$. On the other hand, $\mathbb{Z}/48$ is
detected by ${}_{-1}L_{3}(\mathbb{Z}^{\prime})$.

\textbf{2. }Note that the groups $_{-1}L_{i}(\mathbb{Z}^{\prime})$ for
$i\equiv\pm3\;(\mathrm{\operatorname{mod}\,}8)$ reveal that, in contrast to
the orthogonal situation, the map $B\mathrm{Sp}\mathbb{Z}^{\prime}{}_{\#}%
^{+}\rightarrow B\mathrm{Sp}\mathbb{F}_{3}{}_{\#}^{+}$ above does not admit a section.
\end{Remarks}

\begin{Void}
\textbf{Information on the integral orthogonal group conjecture.}
\end{Void}

The results tabulated above provide a certain amount of support for the
conjectures of Section \ref{background}. To present this information, we
consider the commuting diagram{\tiny
\[%
\begin{array}
[c]{ccccccccccc}%
\text{`}\operatorname{Im}J\text{'} & \rightarrow &  B\Sigma_{\infty}^{+} &
\rightarrow &  B{}_{1}O(\mathbb{Z})^{+} & \rightarrow &  B{}_{1}%
O(\mathbb{Z}^{\prime})^{+} & \rightarrow &  BO\times BO & \overset
{\mathrm{pr}_{1}}{\longrightarrow} & BO\\
&  &  &  & \downarrow^{H} &  & \downarrow^{H} &  &  & \searrow_{\mathrm{sum}%
}^{\mathrm{direct}} & \\
&  &  &  & B\mathrm{GL}(\mathbb{Z})^{+} & \rightarrow &  B\mathrm{GL}%
(\mathbb{Z}^{\prime})^{+} & \rightarrow &  B\mathrm{GL}\mathbb{R}^{\delta+} &
\rightarrow &  BO
\end{array}
\]
}and apply ($2$-primary) input from the following sources. First, the homotopy
groups of $B\Sigma_{\infty}^{+}$ correspond to stable homotopy groups of
spheres by results of \cite{BarrattPriddy}, \cite{Segal: CatsCohom}. In turn,
the map from these groups to the $K$-theory of the integers (presented in
\cite{W: CR}) is studied in \cite{Quillen: to Milnor}. The vanishing of this
map on the image of the $J$-homomorphism in certain dimensions is shown in
\cite{Waldhausen: 1982}. The $KO$-theory \emph{degree map} (induced by group
inclusion) from stable homotopy groups of spheres to the homotopy groups of
$BO$ (known after Bott) is the given by the composite of the upper horizontal
maps, and is analysed in \cite{Adams: J}; it factors through the groups
$_{1}L_{i}(\mathbb{Z}^{\prime})$ computed in Theorem \ref{Thm B in 6} above.
For $i>1$, the map from $B\mathrm{GL}(\mathbb{Z})^{+}$ to $B\mathrm{GL}%
(\mathbb{Z}^{\prime})^{+}$ is known from the localization sequence of
\cite{Soule: Chow} to induce an isomorphism of homotopy groups. Finally,
\cite{S} compares $B\mathrm{GL}\mathbb{R}^{\delta+}$ with $B\mathrm{GL}%
\mathbb{R}\simeq BO$.

Applying $\pi_{i}$ ($i>1$) gives, modulo odd torsion:

\bigskip

\noindent$i\equiv0\ (\mathrm{\operatorname{mod}\,}8).$%
\[%
\begin{array}
[c]{ccccccccccc}%
\mathbb{Z}/2 & \rightarrowtail & ? & \rightarrow & ? & \rightarrow &
\mathbb{Z\oplus Z}/2 & \rightarrow & \mathbb{Z\oplus Z} & \rightarrow &
\mathbb{Z}\\
&  &  &  & \downarrow &  & \downarrow &  &  & \searrow & \\
&  &  &  & 0 & \overset{}{\longrightarrow} & 0 & \rightarrow & 0 & \rightarrow
& \mathbb{Z}%
\end{array}
\]

\bigskip

\noindent$i\equiv1\ (\mathrm{\operatorname{mod}\,}8).$%
\[%
\begin{array}
[c]{ccccccccccc}%
\mathbb{Z}/2 & \rightarrowtail & ? & \rightarrow & ? & \rightarrow &
(\mathbb{Z}/2)^{3} & \rightarrow & \mathbb{Z}/2\mathbb{\oplus Z}/2 &
\rightarrow & \mathbb{Z}/2\\
&  &  &  & \downarrow &  & \downarrow &  &  & \searrow & \\
&  &  &  & \mathbb{Z\oplus Z}/2 & \overset{\mathrm{id}}{\longrightarrow} &
\mathbb{Z\oplus Z}/2 & \rightarrow & \mathbb{Z}/2 & \rightarrow & \mathbb{Z}/2
\end{array}
\]
with the degree map surjective, and composite of top left two arrows followed
by vertical arrow zero.

\bigskip

\noindent$i\equiv2\ (\mathrm{\operatorname{mod}\,}8).$%
\[%
\begin{array}
[c]{ccccccccccc}%
0 & \rightarrow & ? & \rightarrow & ? & \rightarrow & \mathbb{Z}%
/2\mathbb{\oplus Z}/2 & \rightarrow & \mathbb{Z}/2\mathbb{\oplus Z}/2 &
\rightarrow & \mathbb{Z}/2\\
&  &  &  & \downarrow &  & \downarrow &  &  & \searrow & \\
&  &  &  & \mathbb{Z}/2 & \overset{\mathrm{id}}{\longrightarrow} &
\mathbb{Z}/2 & \rightarrow & \mathbb{Z}/2 & \rightarrow & \mathbb{Z}/2
\end{array}
\]
with degree map surjective.

\bigskip

\noindent$i\equiv3\ (\mathrm{\operatorname{mod}\,}8).$%
\[%
\begin{array}
[c]{ccccccccccc}%
\mathbb{Z}/8 & \rightarrowtail & ? & \rightarrow & ? & \rightarrow &
\mathbb{Z}/8 & \rightarrow & 0 & \rightarrow & 0\\
&  &  &  & \downarrow &  & \downarrow &  &  & \searrow & \\
&  &  &  & \mathbb{Z}/16 & \overset{\mathrm{id}}{\longrightarrow} &
\mathbb{Z}/16 & \rightarrow & \mathbb{Q}_{2}/\mathbb{Z}_{2} & \rightarrow & 0
\end{array}
\]
with composite of top three leftmost arrows followed by vertical arrow an injection.

\bigskip

\noindent$i\equiv4\ (\mathrm{\operatorname{mod}\,}8).$%
\[%
\begin{array}
[c]{ccccccccccc}%
0 & \rightarrow & ? & \rightarrow & ? & \rightarrow & \mathbb{Z} & \rightarrow
& \mathbb{Z\oplus Z} & \rightarrow & \mathbb{Z}\\
&  &  &  & \downarrow &  & \downarrow &  &  & \searrow & \\
&  &  &  & 0 & \overset{}{\longrightarrow} & 0 & \rightarrow & 0 & \rightarrow
& \mathbb{Z}%
\end{array}
\]

\bigskip

\noindent$i\equiv5\ (\mathrm{\operatorname{mod}\,}8).$%
\[%
\begin{array}
[c]{ccccccccccc}%
0 & \rightarrow & ? & \rightarrow & ? & \rightarrow & 0 & \rightarrow & 0 &
\rightarrow & 0\\
&  &  &  & \downarrow &  & \downarrow &  &  & \searrow & \\
&  &  &  & \mathbb{Z} & \overset{\mathrm{id}}{\longrightarrow} & \mathbb{Z} &
\rightarrow & 0 & \rightarrow & 0
\end{array}
\]

\bigskip

\noindent$i\equiv6\ (\mathrm{\operatorname{mod}\,}8).$%
\[%
\begin{array}
[c]{ccccccccccc}%
0 & \rightarrow & ? & \rightarrow & ? & \rightarrow & 0 & \rightarrow & 0 &
\rightarrow & 0\\
&  &  &  & \downarrow &  & \downarrow &  &  & \searrow & \\
&  &  &  & 0 & \overset{}{\longrightarrow} & 0 & \rightarrow & 0 & \rightarrow
& 0
\end{array}
\]

\bigskip

\noindent$i\equiv7\ (\mathrm{\operatorname{mod}\,}8),$ $i+1=(\mathrm{odd}%
)\cdot2^{t}.$%
\[%
\begin{array}
[c]{ccccccccccc}%
\mathbb{Z}/2^{t+1} & \rightarrowtail & ? & \rightarrow & ? & \rightarrow &
\mathbb{Z}/2^{t+1} & \rightarrow & 0 & \rightarrow & 0\\
&  &  &  & \downarrow &  & \downarrow &  &  & \searrow & \\
&  &  &  & \mathbb{Z}/2^{t+1} & \overset{\mathrm{id}}{\longrightarrow} &
\mathbb{Z}/2^{t+1} & \rightarrow & \mathbb{Q}_{2}/\mathbb{Z}_{2} & \rightarrow
& 0
\end{array}
\]
with composite of top three leftmost arrows followed by vertical arrow an isomorphism.

\section{\label{Homotopy fixed points}Homotopy fixed points}

In response to a query of B. Kahn, we now investigate the homotopy fixed point
set of $\mathbb{Z}/2$ acting on $\mathcal{K}(\mathbb{Z}^{\prime})$, at least
after $2$-adic completion. We recall that, for a ring $\Lambda$,
$\mathcal{K}(\Lambda)$ is described \emph{functorially }as the loop space
$\Omega B\mathrm{GL}(S\Lambda)^{+}$. Therefore our main interest in the
description below is in the case $A=S\Lambda$ (see Appendix A), and in
particular $\Lambda=\mathbb{Z}^{\prime}$.

In order to describe the $\mathbb{Z}/2$ action with respect to a ring $A$
admitting an antiinvolution $x\mapsto\check{x}$, we recall that, for
$\varepsilon=\pm1$ fixed by the involution, ${}_{\varepsilon}O_{n,n}(A)$ as
defined above is in effect the fixed subgroup $\mathrm{GL}_{2n}(A)^{{}%
_{\varepsilon}\mathbb{Z}/2}$ of the action denoted ${}_{\varepsilon}%
\mathbb{Z}/2$ on $\mathrm{GL}_{2n}(A)$. (Our notation omits $\varepsilon$ when
it is safe to do so.) The nontrivial element of ${}_{\varepsilon}\mathbb{Z}/2$
sends an invertible matrix to the $\varepsilon$-hyperbolic adjoint of its
inverse. Explicitly,
\[
M=\left(
\begin{array}
[c]{cc}%
a & b\\
c & d
\end{array}
\right)  \in\mathrm{GL}_{2n}(A)
\]
is sent to the inverse of
\[
{}_{\varepsilon}J_{n}\cdot{}^{\mathrm{t}}\check{M}\cdot{}_{\varepsilon}%
J_{n}^{-1}=\left(
\begin{array}
[c]{cc}%
^{\mathrm{t}}\check{d} & \varepsilon\,{}^{\mathrm{t}}\check{b}\\
\varepsilon\,{}^{\mathrm{t}}\check{c} & ^{\mathrm{t}}\check{a}%
\end{array}
\right)
\]
where
\[
{}_{\varepsilon}J_{n}=\left(
\begin{array}
[c]{cc}%
0 & \varepsilon I_{n}\\
I_{n} & 0
\end{array}
\right)  \text{.}%
\]
We denote by ${}_{\varepsilon}O(A)$ the direct limit of the ${}_{\varepsilon
}O_{n,n}(A)$ with respect to the obvious inclusions within each of the four
component blocks. If $\Lambda=\mathbb{R}$ or $\mathbb{C}$, we take for
$A=S\Lambda$ the topological suspension (as described in Appendix A).

As in \cite{Faj}, one observes that the ${}_{\varepsilon}\mathbb{Z}/2$ action
is compatible with the process of stabilization, passage to classifying spaces
and the plus-construction. Moreover, it behaves well with respect to
suspension, and passing to loop spaces commutes with taking fixed point sets.
Thus the fixed point set of the space $\mathcal{K}(\Lambda)=\Omega
(B\mathrm{GL}(S\Lambda)^{+})$ is just%
\[
{}_{\varepsilon}\mathcal{L}(\Lambda)=\Omega(B{}_{\varepsilon}O(S\Lambda
)^{+})\text{,}%
\]
a fact that we record.

\begin{Lemma}
For $\Lambda$ as above, and $\varepsilon=\pm1$, the fixed point set of the
${}_{\varepsilon}\mathbb{Z}/2$ action on $\mathcal{K}(\Lambda)$ is
\[
\mathcal{K}(\Lambda){}^{{}_{\varepsilon}\mathbb{Z}/2}\simeq{}_{\varepsilon
}\mathcal{L}(\Lambda)\text{.}%
\]
\hfill$\Box$
\end{Lemma}

We now turn our attention to the homotopy fixed point set
\[
\mathcal{K}(\Lambda){}^{h({}_{\varepsilon}\mathbb{Z}/2)}:=\mathrm{map}%
_{{}_{\varepsilon}\mathbb{Z}/2}(E\mathbb{Z}/2,\,\mathcal{K}(\Lambda))\text{,}%
\]
the space of maps equivariant under the ${}_{\varepsilon}\mathbb{Z}/2$ action,
where $E\mathbb{Z}/2$ is a contractible free $\mathbb{Z}/2$-space (usually
taken to be $S^{\infty}$ with antipodal action), and $\Lambda=\mathbb{Z}%
^{\prime}$ and its related rings, such as $\mathbb{R}$, $\mathbb{C}$ or
$\mathbb{F}_{3}$.

We shall prove the following theorem.

\begin{theorem}
\label{Thm C in 7}For $\varepsilon=\pm1$, the natural map%
\[
{}_{\varepsilon}\mathcal{L}(\mathbb{Z}^{\prime})_{\#}\longrightarrow
\lbrack\mathcal{K}(\mathbb{Z}^{\prime})^{h({}_{\varepsilon}\mathbb{Z}%
/2)}]_{\#}%
\]
is a homotopy equivalence.
\end{theorem}

Since the spaces ${}_{\varepsilon}\mathcal{L}(\mathbb{Z}^{\prime})_{\#}$ and
$[\mathcal{K}(\mathbb{Z}^{\prime})^{h({}_{\varepsilon}\mathbb{Z}/2)}]_{\#}$
are both homotopy pull-backs, the proof largely reduces to three pairs of
verifications (Lemmas \ref{hZ/2 for C}, \ref{hZ/2 for F_3}, \ref{hZ/2 for R}
below) which may be found to some extent in the literature (\cite{Faj},
\cite{HopMahSad}, \cite{K: descent}, \cite{Kobal: thesis}, \cite{Kobal}).
Because the spaces of the theorem are not connected, the group of connected
components needs special attention here.

\begin{Lemma}
\label{hZ/2 for C}For $\varepsilon=\pm1$, the natural map
\[
{}_{\varepsilon}\mathcal{L}(\mathbb{C})\longrightarrow\mathcal{K}%
(\mathbb{C)}^{h({}_{\varepsilon}\mathbb{Z}/2)}%
\]
is a homotopy equivalence.
\end{Lemma}

\noindent\textbf{Proof. }The lemma is true without the need of $2$-adic
completion. The case $\varepsilon=1$ has a simple analysis, by means of the
map of connected-component fibrations $X^{0}\rightarrow X\rightarrow\pi
_{0}(X)$:
\[%
\begin{array}
[c]{ccccc}%
\Omega^{0}B{}_{1}O(S\mathbb{C})=B{}_{1}O(\mathbb{C}) & \longrightarrow &
{}_{1}\mathcal{L}(\mathbb{C}) & \longrightarrow & {}_{1}L_{0}(\mathbb{C})\\
\downarrow &  & \downarrow &  & \downarrow\\
\Omega^{0}B\mathrm{GL}(S\mathbb{C})=B\mathrm{GL}(\mathbb{C}) & \longrightarrow
& \mathcal{K}(\mathbb{C}) & \longrightarrow &  K_{0}(\mathbb{C})
\end{array}
\]
that we claim represents the inclusion of the homotopy fixed point sets.
Certainly this is true for the discrete base $K_{0}(\mathbb{C})$, since the
action of $\mathbb{Z}/2$ is trivial and the map $\mathbb{Z=\,}{}_{1}%
L_{0}(\mathbb{C})\rightarrow K_{0}(\mathbb{C})=\mathbb{Z}$ is the identity. By
polar decomposition of matrices, the required result for the fiber follows
from a classical fact in algebraic topology: the map%
\[
BO\longrightarrow BU^{h({}_{1}\mathbb{Z}/2)}%
\]
is a homotopy equivalence \cite{HopMahSad}, \cite{Faj}, \cite{K: descent}.

For $\varepsilon=-1$, we cannot repeat the method above because the map
${}_{-1}L_{0}(\mathbb{C})\rightarrow K_{0}(\mathbb{C})$ fails to be an
isomorphism. Instead we use the general argument of \cite{K: descent} for
Banach algebras. Specifically, we use the fact that, because $A=S\mathbb{C}$
is a $C^{\ast}$-algebra, the inclusion of $U_{2n}(A)$ in $\mathrm{GL}_{2n}(A)$
is a $\mathbb{Z}/2$-homotopy equivalence, where
\[
U_{2n}(A)=\left\{  M=\left(
\begin{array}
[c]{cc}%
a & b\\
c & d
\end{array}
\right)  \mid M^{-1}=\left(
\begin{array}
[c]{cc}%
^{\mathrm{t}}\bar{a} & ^{\mathrm{t}}\bar{c}\\
^{\mathrm{t}}\bar{b} & ^{\mathrm{t}}\bar{d}%
\end{array}
\right)  \right\}
\]
and the involution on $\mathrm{GL}_{2n}(A)$ is
\[
M\longmapsto\left(
\begin{array}
[c]{cc}%
^{\mathrm{t}}d & -\,{}^{\mathrm{t}}b\\
-\,{}^{\mathrm{t}}c & ^{\mathrm{t}}a
\end{array}
\right)  ^{-1}\text{.}%
\]
So, when $M\in U_{2n}(A)$, it is sent to
\[
\left(
\begin{array}
[c]{cc}%
\bar{d} & -\bar{c}\\
-\bar{b} & \bar{a}%
\end{array}
\right)  ={}_{-1}J_{n}\bar{M}{}_{-1}J_{n}^{-1}\text{.}%
\]
Therefore $U_{2n}(A)$ is in turn $\mathbb{Z}/2$-homotopy equivalent to
$\mathrm{GL}_{2n}(A)$ equip\-ped with the involution $M\mapsto{}_{-1}J_{n}%
\bar{M}{}_{-1}J_{n}^{-1}$. If we reinterpret $\mathrm{GL}_{2n}(S\mathbb{C})$
as $\mathrm{GL}_{n}(S\mathbb{H\otimes}_{\mathbb{R}}\mathbb{C})$, with basis
afforded by
\[
I_{2},\ \left(
\begin{array}
[c]{cc}%
0 & 1\\
-1 & 0
\end{array}
\right)  ,\ \left(
\begin{array}
[c]{cc}%
i & 0\\
0 & -i
\end{array}
\right)  ,\ \left(
\begin{array}
[c]{cc}%
0 & -i\\
-i & 0
\end{array}
\right)
\]
then the involution now acts as complex conjugation. Therefore, by the main
theorem of\textsf{ }\cite{K: descent}, there is a homotopy equivalence between
$\mathcal{K}(S\mathbb{H})$ and $\mathcal{K}(S\mathbb{C})^{h({}_{-1}%
\mathbb{Z}/2)}$, and hence, on passing to loop spaces, between $\mathcal{K}%
(\mathbb{H})={}_{-1}\mathcal{L}(\mathbb{C})$ and $\mathcal{K}(\mathbb{C}%
)^{h({}_{-1}\mathbb{Z}/2)}$.\hfill$\Box$

$\smallskip$

\begin{Lemma}
\label{hZ/2 for F_3}For $\varepsilon=\pm1$ and $q$ odd, the natural map
\[
{}_{\varepsilon}\mathcal{L}(\mathbb{F}_{q})\longrightarrow\mathcal{K}%
(\mathbb{F}_{q})^{h({}_{\varepsilon}\mathbb{Z}/2)}%
\]
is a homotopy equivalence.
\end{Lemma}

\noindent\textbf{Proof. }This result also is true without the need of $2$-adic completion.

From Theorem \ref{Friedlander fibn} and fibration (\ref{QF fibrations})
above\textsf{, }there is a map of homotopy fibrations%
\[%
\begin{array}
[c]{ccccc}%
{}_{\varepsilon}\mathcal{L}(\mathbb{F}_{q}) & \longrightarrow & {}%
_{\varepsilon}\mathcal{L}(\mathbb{C}) & \overset{\psi^{q}-1}{\longrightarrow}%
& {}_{\varepsilon}\mathcal{L}(\mathbb{C})\\
\downarrow &  & \downarrow &  & \downarrow\\
\mathcal{K}(\mathbb{F}_{q})^{h({}_{\varepsilon}\mathbb{Z}/2)} &
\longrightarrow & \mathcal{K}(\mathbb{C})^{h({}_{\varepsilon}\mathbb{Z}/2)} &
\overset{\psi^{q}-1}{\longrightarrow} & \mathcal{K}(\mathbb{C})^{h({}%
_{\varepsilon}\mathbb{Z}/2)}%
\end{array}
\]
So the result follows from our previous discussion of the complex case. See
also \cite{Kobal: thesis}.\hfill$\Box\smallskip$

To complete the proof of Theorem \ref{Thm C in 7}, we now turn to
consideration of the homotopy fixed point set for $\mathbb{R}$. This
verification is the most tricky one and really needs the $2$-adic completion.
Again, although this result may be found in the literature (at least, in
positive dimensions), we present an alternative viewpoint using Fredholm
operators that we feel is illuminating. The result that we require is the
following. (Since we are considering $A=S\mathbb{R}$, $\mathcal{K}%
(\mathbb{R})=\Omega B\mathrm{GL}(A)^{+}$ simplifies to $\mathrm{GL}(A)$;
similarly for ${}_{\varepsilon}\mathcal{L}(\mathbb{R})$.)

\begin{Lemma}
\label{hZ/2 for R}The natural map
\[
{}_{\varepsilon}\mathcal{L}(\mathbb{R})={}_{\varepsilon}O(A)\longrightarrow
\mathrm{GL}(A)^{h\mathbb{Z}/2}=\mathcal{K}(\mathbb{R})^{h\mathbb{Z}/2}%
\]
becomes a homotopy equivalence after $2$-adic completion.
\end{Lemma}

We now prove this lemma for $\varepsilon=1$, and then use elements of this
argument when we later establish the case $\varepsilon=-1$.

\begin{Void}
\textbf{The case} $\varepsilon=1$\textbf{.}
\end{Void}

This case calls for analysis of the space ${}_{1}\mathcal{L}(\mathbb{R}%
)={}_{1}O(A)$, where $A=S\mathbb{R}$. As above, we regard ${}_{1}O(A)$ as the
direct limit of the ${}_{1}O_{n,n}(A)$. Using polar decomposition on matrices
over the $C^{\ast}$-algebra $A$, we reinterpret ${}_{1}O_{n,n}(A)$ (up to
homotopy) as the group of fixed points in $O_{2n}(A)$ under conjugation by the
matrix $\left(
\begin{array}
[c]{cc}%
I_{n} & 0\\
0 & -I_{n}%
\end{array}
\right)  $. In other words, we have a commutative diagram of $\mathbb{Z}%
/2$-homotopy equivalence inclusions%
\[%
\begin{array}
[c]{ccc}%
{}_{1}O_{n,n}(A) & \hookrightarrow & \mathrm{GL}_{2n}(A)\\
\uparrow &  & \uparrow\\
O_{n}(A)\times O_{n}(A) & \hookrightarrow &  O_{2n}(A)\\
\downarrow &  & \downarrow\\
\mathrm{GL}_{n}(A)\times\mathrm{GL}_{n}(A) & \hookrightarrow & \mathrm{GL}%
_{2n}(A)
\end{array}
\]
Since by Kuiper's theorem the inclusion $\mathrm{GL}_{n}(A)\hookrightarrow
\mathrm{GL}_{n+1}(A)$ is a homotopy equivalence, we just have to investigate
the map%
\[
\mathrm{GL}_{1}(A)\times\mathrm{GL}_{1}(A)\longrightarrow\mathrm{GL}%
_{2}(A)^{h\mathbb{Z}/2}%
\]
or equivalently (as in Appendix A)
\[
\chi:\mathcal{F}(H)\times\mathcal{F}(H)\longrightarrow\mathcal{F}(H\oplus
H)^{h\mathbb{Z}/2}\text{.}%
\]
Here $\chi$ sends $(D_{1},\,D_{2})$ to the constant map $E\mathbb{Z}%
/2\rightarrow\mathcal{F}(H\oplus H)$ that takes value $D_{1}\oplus D_{2}$.
Recall from the Atiyah-J\"{a}nich theorem (see for example \cite{Atiyah: book}
or \cite{Karoubi: Espaces classifiants}) that the space $\mathcal{F}(H)$ of
Fredholm operators in a (real or complex) Hilbert space $H$ serves as a
classifying space for the $K$-theory $K(X)$ of paracompact spaces $X$.

The statement that we want to prove for $\varepsilon=1$ therefore amounts to
the assertion that%
\[
\mathcal{K}(\mathbb{R})\times\mathcal{K}(\mathbb{R})\longrightarrow
\mathcal{K}(M_{2}(\mathbb{R}))^{h\mathbb{Z}/2}%
\]
becomes a $2$-adic equivalence, where the $\mathbb{Z}/2$-action on
$M_{2}(S\mathbb{R})$ is given by conjugation with respect to $\left(
\begin{array}
[c]{cc}%
1 & 0\\
0 & -1
\end{array}
\right)  $. For later use, we also need to prove that
\[
\mathcal{K}(\mathbb{C})\times\mathcal{K}(\mathbb{C})\longrightarrow
\mathcal{K}(M_{2}(\mathbb{C}))^{h\mathbb{Z}/2}%
\]
becomes a $2$-adic equivalence. Therefore in what follows, $H$ represents a
real or complex Hilbert space. Also, we write $\mathrm{map}$ for the
nonequivariant mapping space, and $\mathrm{map}_{\ast}$ for the subspace of
basepoint-preserving maps.

\begin{Lemma}
\label{2-adic BO}The function space $\mathrm{map}_{\ast}(B\mathbb{Z}%
/2,\,\mathcal{F}(H))$ is a version $\mathcal{F}(H)_{\#}$ of the $2$-adic
completion of $\mathcal{F}(H)$ (including the completion of $\pi_{0}$). With
this version, the canonical map $\mathcal{F}(H)\rightarrow\mathcal{F}(H)_{\#}$
sends an element $x$ of $K(X)=[X,\,\mathcal{F}(H)]$ to the map from $X\times
B\mathbb{Z}/2$ to $\mathcal{F}(H)$ associated to the element $x(L-1)$ of the
relative $K$-group $K(X\times B\mathbb{Z}/2,X)$, where $L$ is the canonical
line bundle over $B\mathbb{Z}/2$.
\end{Lemma}

\noindent\textbf{Proof.} Let $X$ be a sphere. We have to show that the map
$\mathcal{F}(H)\rightarrow\mathrm{map}_{\ast}(B\mathbb{Z}/2,\,\mathcal{F}(H))$
described above induces an isomorphism%
\[
\lbrack X,\,\mathcal{F}(H)]_{\#}\rightarrow\lbrack X,\,\mathrm{map}_{\ast
}(B\mathbb{Z}/2,\,\mathcal{F}(H))]
\]
where $[X,\,\mathcal{F}(H)]_{\#}$ is the $2$-adic completion of the group
$[X,\,\mathcal{F}(H)]$. The function space $\mathrm{map}(X,\,\mathrm{map}%
_{\ast}(B\mathbb{Z}/2,\,\mathcal{F}(H)))$ may be identified with the function
space of maps $X\times B\mathbb{Z}/2\rightarrow\mathcal{F}(H)$ that send
$X\times\{\ast\}$ to the basepoint $1_{H}$ of $\mathcal{F}(H)$. In other
words, we must calculate the relative group $K(X\times\mathbb{R}P^{\infty
},\,X)$, and investigate the map%
\[
K(X)\longrightarrow K(X\times\mathbb{R}P^{\infty},\,X)
\]
given by the correspondence $x\mapsto x(L-1)$ as in the statement of the lemma.

Since the groups $K(X\times\mathbb{R}P^{n},\,X)$ are finite for $n$ even (by,
for example, the Atiyah-Hirzebruch spectral sequence, or the exact sequences
below), the projective limit $\underleftarrow{\mathrm{\lim}}K(X\times
\mathbb{R}P^{n},\,X)$ has vanishing $\underleftarrow{\mathrm{\lim}}^{1}$ term
and is therefore isomorphic to $K(X\times\mathbb{R}P^{\infty},\,X)$.

Now the group $K(X\times\mathbb{R}P^{n},\,X)$ has been described concretely in
\cite{Kar}p.249, Theorem 6.40 (see also, more recently, \cite{K:
equivariant}). It is the middle term of an exact sequence{\small
\[
K(\mathcal{E}^{n,0}(X))\rightarrow K(\mathcal{E}(X))\overset{\alpha
}{\rightarrow}K(X\times\mathbb{R}P^{n},\,X)\rightarrow K^{1}(\mathcal{E}%
^{n,0}(X))\overset{\Phi_{n}}{\longrightarrow}K^{1}(\mathcal{E}(X))\text{,}%
\]
}where $\mathcal{E}(X)$ denotes the category of (real or complex) vector
bundles on $X$, and $\mathcal{E}^{n,0}(X)$ is the category of Clifford bundles
(those with an action of the usual Clifford algebra $C^{n,0}$). Note that
$\alpha$ is the cup-product by $(L-1)$, where $L$ is the canonical line bundle
over $\mathbb{R}P^{n}$. Since we are taking the projective limit, we may
assume that $n=8k$. Then the Clifford algebra $C^{n,0}$ is isomorphic to the
matrix algebra $M_{2^{4k}}(\mathbb{R})$, the category $\mathcal{E}^{n,0}(X)$
is equivalent to $\mathcal{E}(X)$, and the `restriction of scalars' functor
$\mathcal{E}^{n,0}(X)\longrightarrow\mathcal{E}(X)$ simply sends the class of
a bundle $E$ to that of its multiple $2^{4k}E$. On the other hand, since $X$
is a sphere, $K^{1}(\mathcal{E}^{n,0}(X))=K^{1}(X)$ is $\mathbb{Z}$ or $0$ (in
the complex case), else $\mathbb{Z}$, $0$ or $\mathbb{Z}/2$ (in the real case).

From the diagram{\tiny
\[%
\begin{array}
[c]{ccccccccc}%
K(\mathcal{E}^{n+8,0}(X)) & \rightarrow &  K(\mathcal{E}(X)) & \rightarrow &
K(X\times\mathbb{R}P^{n+8},\,X) & \rightarrow &  K^{1}(\mathcal{E}%
^{n+8,0}(X)) & \overset{\Phi_{n+8}}{\longrightarrow} & K^{1}(\mathcal{E}(X))\\
\downarrow^{\cdot16} &  & \downarrow^{\mathrm{id}} &  & \downarrow &  &
\downarrow^{\cdot16} &  & \downarrow^{\mathrm{id}}\\
K(\mathcal{E}^{n,0}(X)) & \rightarrow &  K(\mathcal{E}(X)) & \rightarrow &
K(X\times\mathbb{R}P^{n},\,X) & \rightarrow &  K^{1}(\mathcal{E}^{n,0}(X)) &
\overset{\Phi_{n}}{\longrightarrow} & K^{1}(\mathcal{E}(X))
\end{array}
\]
}we see that the induced map $\mathrm{Ker}\Phi_{n+8}\rightarrow\mathrm{Ker}%
\Phi_{n}$ is always zero. Therefore, $K(X\times\mathbb{R}P^{\infty
},\,X)=\underleftarrow{\mathrm{\lim}}K(X\times\mathbb{R}P^{n},\,X)$ is
isomorphic to
\[
\underleftarrow{\mathrm{\lim}}K(X)/2^{4k}K(X)=K(X)_{\#}.
\]
\hfill$\Box\smallskip$

\begin{Corollary}
\label{split fibn}The above splitting
\[
K(X\times\mathbb{R}P^{\infty})=K(X)\oplus K(X\times\mathbb{R}P^{\infty},\,X)
\]
corresponds to the splitting of the function space $\mathrm{map}%
(B\mathbb{Z}/2,\,\mathcal{F}(H))$ as the product
\[
\theta:\mathrm{map}(B\mathbb{Z}/2,\,\mathcal{F}(H))\overset{\mathrm{\simeq}%
}{\longrightarrow}\mathcal{F}(H)\times\mathrm{map}_{\ast}(B\mathbb{Z}%
/2,\,\mathcal{F}(H))=\mathcal{F}(H)\times\mathcal{F}(H)_{\#}\text{.}%
\]
\end{Corollary}

\noindent\textbf{Proof. }In view of Appendix A, the split fibration%
\[
\mathrm{map}_{\ast}(B\mathbb{Z}/2,\,\mathcal{F}(H))\longrightarrow
\mathrm{map}(B\mathbb{Z}/2,\,\mathcal{F}(H))\overset{\mathrm{eval}%
}{\longrightarrow}\mathcal{F}(H)
\]
can be interpreted as the split fibration of topological groups (hence a
product fibration as spaces)%
\[
\mathrm{map}_{\ast}(B\mathbb{Z}/2,\,\mathrm{GL}_{1}(A))\longrightarrow
\mathrm{map}(B\mathbb{Z}/2,\,\mathrm{GL}_{1}(A))\overset{\mathrm{eval}%
}{\longrightarrow}\mathrm{GL}_{1}(A)
\]
where $\mathrm{GL}_{1}(A)$ is the group of invertible elements in the Calkin
algebra $S\mathbb{R}$ or $S\mathbb{C}$, and the projection map is evaluation
at the basepoint. One has to recall also that the group structure on
$\mathrm{GL}_{1}(A)$ induces the group structure on $K(X)=[X,\,\mathrm{GL}%
_{1}(A)]$.\hfill$\Box\smallskip$

This result combines with the classical homotopy equivalence%
\[
\mathcal{F}(H)\simeq\mathcal{F}(H\oplus H)
\]
of Kuiper \cite{Kuiper} and Palais \cite{Palais} to yield the following.

\begin{Corollary}
\label{F(H+H)}The above splitting induces a splitting
\[
\Psi:\mathcal{F}(H\oplus H)^{h\mathbb{Z}/2}\overset{\simeq}{\longrightarrow
}\mathcal{F}(H)\times\mathrm{map}_{\ast}(B\mathbb{Z}/2,\,\mathcal{F}%
(H))=\mathcal{F}(H)\times\mathcal{F}(H)_{\#}.
\]
\hfill$\Box$
\end{Corollary}

We now have the composition%
\begin{align*}
\mathcal{F}(H)\times\mathcal{F}(H)  &  \overset{\chi}{\longrightarrow
}\mathcal{F}(H\oplus H)^{h\mathbb{Z}/2}\\
&  \overset{\simeq}{\longrightarrow}\mathcal{F}(H)\times\mathrm{map}_{\ast
}(B\mathbb{Z}/2,\,\mathcal{F}(H))=\mathcal{F}(H)\times\mathcal{F}(H)_{\#}%
\end{align*}
which we seek to show becomes an equivalence on $2$-adic completion. This
composite is analysed in terms of its effect on groups of homotopy classes of
maps from a compact space $X$ into the function spaces.

For this, we recall that another description of $K(X)=[X,\,\mathcal{F}(H)]$
consists of homotopy classes $d(E,F,D)$ of triples $(E,F,D)$ where $E$ and $F$
are infinite-dimensional Hilbert bundles and $D:E\rightarrow F$ is a bundle
morphism such that $D_{x}:E_{x}\rightarrow F_{x}$ is a Fredholm map for each
$x\in X$ (see for instance \cite{Kar}\S\ II.2 for the general framework). As a
matter of fact, Kuiper's theorem on the contractibility of $\mathrm{Aut}(H)$
implies that, without loss of generality, one can take $E,F$ to be trivial
Hilbert bundles. In this way $D$ defines a map $X\rightarrow\mathcal{F}(H)$
which is unique up to homotopy. One advantage of this presentation is the
following: for a finite-dimensional vector bundle $G$ defining a class in
$K(X)$, the cup-product of $x=d(E,F,D)$ by $G$ is simply $d(E\otimes
G,\,F\otimes G,\,D\otimes\mathrm{id}_{G})$. Note also that the isomorphism
$K(X)\rightarrow\lbrack X,\,\mathcal{F}(H)]$ is just induced by the
correspondence associating to a vector bundle $G$ over $X$ the class
$d(E,F,D)$ where
\[
E=F=\ell^{2}(G)=G\oplus\cdots\oplus G\oplus\cdots
\]
and
\[
D(x_{0},x_{1},x_{2},\ldots)=(x_{1},x_{2},\ldots)\text{.}%
\]
Finally, for the group structure in $K(X)$, we have%
\[
d(E,F,D)+d(E^{\prime},F^{\prime},D^{\prime})=d(E\oplus E^{\prime},F\oplus
F^{\prime},D\oplus D^{\prime})\text{.}%
\]

That $\chi$ becomes an equivalence on $2$-completion is immediate from the
explicit description of the following lemma.
%
%
%
%
%

\begin{Lemma}
\label{(E,F)}The natural map%
\nopagebreak
\[
\Psi\circ\chi:\mathcal{F}(H)\times\mathcal{F}(H)\longrightarrow\mathcal{F}%
(H\oplus H)^{h\mathbb{Z}/2}\overset{\simeq}{\longrightarrow}\mathcal{F}%
(H)\times\mathrm{map}_{\ast}(B\mathbb{Z}/2,\,\mathcal{F}(H))
\]%
\nopagebreak
is induced by the correspondence (seen from the point of view of $K$-theory,
as in Lemma \ref{split fibn}) that associates to a pair of vector bundles
$(E_{1},E_{2})$ the virtual bundles $E_{1}\oplus E_{2}$ and $E_{2}%
\otimes(L-\mathbf{1})$, where $L$ is the Hopf line bundle and $\mathbf{1}$ is
the trivial line bundle on $B\mathbb{Z}/2=\mathbb{R}P^{\infty}$. Thus, from
Lemma \ref{2-adic BO}, $\Psi\circ\chi$ is a $2$-adic equivalence.
\end{Lemma}

\noindent\textbf{Proof.}\quad We check commutativity of the diagram%
\[%
\begin{array}
[c]{ccccc}%
&  & \mathcal{F}(H) &  & \\
& \swarrow^{\mathrm{in}_{1}}\quad &  & ^{\rho_{1}}\searrow & \\
\mathcal{F}(H)\times\mathcal{F}(H) & \quad\overset{\chi}{\longrightarrow} &
\mathcal{F}(H\oplus H)^{h\mathbb{Z}/2} & \overset{\gamma}{\longleftarrow}\quad
& \mathrm{map}(B\mathbb{Z}/2,\,\mathcal{F}(H))\\
& \nwarrow_{\mathrm{in}_{2}}\quad &  & _{\rho_{L}}\nearrow & \\
&  & \mathcal{F}(H) &  &
\end{array}
\]
and then discuss the effect of composition of the diagram with the canonical
map
\[
\theta:\mathrm{map}(B\mathbb{Z}/2,\,\mathcal{F}(H))\rightarrow\mathcal{F}%
(H)\times\mathrm{map}_{\ast}(B\mathbb{Z}/2,\,\mathcal{F}(H)).
\]
Here $\rho_{1}(E)=E\otimes\mathbf{1}$ and $\rho_{L}(E)=E\otimes L$; we also
write $\chi_{j}=\chi\circ\mathrm{in}_{j}$.

By Corollary \ref{F(H+H)}, the obvious map%
\[
H\longrightarrow H\oplus0\hookrightarrow H\oplus H
\]
induces a homotopy equivalence%
\[
\gamma:\mathrm{map}(B\mathbb{Z}/2,\,\mathcal{F}(H))=\mathcal{F}(H\oplus
0)^{h\mathbb{Z}/2}\longrightarrow\mathcal{F}(H\oplus H)^{h\mathbb{Z}%
/2}\text{,}%
\]
and so $\gamma^{-1}$ is defined up to homotopy.

With this identification, the composition
\[
\mathcal{F}(H)\overset{\chi_{1}}{\longrightarrow}\mathcal{F}(H\oplus
H)^{h\mathbb{Z}/2}\overset{\gamma^{-1}}{\longrightarrow}\mathrm{map}%
(B\mathbb{Z}/2,\,\mathcal{F}(H))
\]
associates to a vector bundle $E_{1}$ over a space $X$ the vector bundle
$\rho_{1}(E_{1})=E_{1}\otimes\mathbf{1}$ on the space $X\times B\mathbb{Z}/2$.

However, the crucial point is to make explicit the image of the `second
component' $E_{2}$, that is, to specify the composition
\[
\mathcal{F}(H)\overset{\chi_{2}}{\longrightarrow}\mathcal{F}(H\oplus
H)^{h\mathbb{Z}/2}\overset{\gamma^{-1}}{\longrightarrow}\mathrm{map}%
(B\mathbb{Z}/2,\,\mathcal{F}(H))
\]
where $\chi_{2}(D)=\mathrm{id}_{H}\oplus D$.

In order to make $\gamma^{-1}\circ\chi_{2}$ explicit, we also translate this
problem geometrically in terms of vector bundles: we start with the class
$x=d(T,T,D)$ where $T$ is the trivial Hilbert bundle $X\times H$, and $D$ is
represented by an endomorphism of $T$ such that the induced map from $X$ to
$\mathrm{End}(H)$ is a family of Fredholm operators. By Kuiper's theorem that
the topological group $\mathrm{Aut}(H)$ is contractible, we construct by
induction on $n$ a map $\sigma:S^{n}\longrightarrow\mathrm{Aut}(H)$ such that
$\sigma(e)=\mathrm{id}_{H}$ and $\sigma(-x)=-\sigma(x)$, where $e$ is the
basepoint of the sphere. We denote again by $\sigma:S^{\infty}\longrightarrow
\mathrm{Aut}(H)$ the map on the infinite sphere (recall that $B\mathbb{Z}%
/2=\mathbb{R}P^{\infty}=S^{\infty}/(\mathbb{Z}/2)$). Now consider the
equivariant homotopy $\lambda(x)\Delta\lambda(x)^{-1}$ in the space of
Fredholm operators on $H\oplus H$, where $\lambda(x)$ and $\Delta$ are given
by the following matrices:%
\[
\lambda(x)=\left(
\begin{array}
[c]{cc}%
\cos\theta & -\sigma(x)\sin\theta\\
\sigma(x)^{-1}\sin\theta & \cos\theta
\end{array}
\right)  \,\quad\text{and\thinspace\quad}\Delta=\left(
\begin{array}
[c]{cc}%
I & 0\\
0 & D
\end{array}
\right)
\]
and $\theta\in\lbrack0,\pi/2]$. At the end of the homotopy we get the family
of Fredholm operators%
\[
\sigma(x)D\sigma(x)^{-1}\oplus\mathrm{id}_{H}\text{.}%
\]
In our geometric interpretation the trick now is to remark that $\sigma(x)$
defines an isomorphism between the trivial Hilbert bundle $T=B\mathbb{Z}%
/2\times H$ and its tensor product $T\otimes L$ with the Hopf bundle $L$ on
the infinite projective space (note that $T\otimes L$ is isomorphic to $T$,
but not canonically). The result for $\gamma^{-1}\circ\chi_{2}(x)$ is then
$d(T\otimes L,\,T\otimes L,\,D\otimes\mathrm{id}_{L})$. In other words, a
vector bundle $E_{2}$ on $X$ goes to $\rho_{L}(E_{2})=$ $E_{2}\otimes L$ on
$X\times B\mathbb{Z}/2$.

Next recall from Corollary \ref{split fibn} that $\theta$ is geometrically
defined by splitting the evaluation map $\mathrm{map}(B\mathbb{Z}%
/2,\,\mathcal{F}(H))\rightarrow\mathrm{map}(\ast,\,\mathcal{F}(H))$. Thus, to
a triple $d(E,F,D)$ as above ($E$ and $F$ being Hilbert bundles on $X\times
B\mathbb{Z}/2$) it associates the pair
\[
(d(E^{\prime},F^{\prime},D^{\prime}),\;d(E,F,D)-d(E^{\prime},F^{\prime
},D^{\prime})),
\]
where $D^{\prime}:E^{\prime}\rightarrow F^{\prime}$ is the pull-back of
$D:E\rightarrow F$ on $X\times B\mathbb{Z}/2$ via the composition map
\[
X\times B\mathbb{Z}/2\rightarrow X\times\ast\rightarrow X\times B\mathbb{Z}%
/2\text{,}%
\]
and $\ast$ is the basepoint of $B\mathbb{Z}/2$. Therefore, according to this
geometric interpretation of $\theta$, we have that $\theta(E_{1}%
\otimes\mathbf{1})=(E_{1},0)$ and $\theta(E_{2}\otimes L)=(E_{2}%
,\,E_{2}\otimes L-E_{2}\otimes\mathbf{1})$. Consequently, on summing we find
that the composition $\Psi\circ\chi=\theta\circ\gamma^{-1}\circ\chi$ sends
$(E_{1},E_{2})$ to $(E_{1}\oplus E_{2},\,E_{2}\otimes(L-\mathbf{1}))$%
.\hfill$\Box$

$\smallskip$

\begin{Void}
\textbf{The case} $\varepsilon=-1$\textbf{.}
\end{Void}

Here we show that ${}_{-1}\mathcal{L}(\mathbb{R})_{\#}\simeq(\mathcal{K}%
(\mathbb{R})^{h({}_{-1}\mathbb{Z}/2)})_{\#}$. Using polar decomposition of
matrices again, we are reduced to showing that
\begin{equation}
\mathcal{K}(M_{2}(\mathbb{R}))_{\#}^{\sigma}\simeq\mathcal{K}(M_{2}%
(\mathbb{R}))_{\#}^{h\sigma} \tag{S}\label{K(M2(R)) fixed points}%
\end{equation}
where $\sigma$ is the involution on $M_{2}(\mathbb{R})$ defined by%
\[
\left(
\begin{array}
[c]{cc}%
a & b\\
c & d
\end{array}
\right)  \longmapsto\left(
\begin{array}
[c]{cc}%
d & -c\\
-b & a
\end{array}
\right)  =\left(
\begin{array}
[c]{cc}%
0 & -1\\
1 & 0
\end{array}
\right)  \left(
\begin{array}
[c]{cc}%
a & b\\
c & d
\end{array}
\right)  \left(
\begin{array}
[c]{cc}%
0 & 1\\
-1 & 0
\end{array}
\right)  \text{.}%
\]
Note that $M_{2}(\mathbb{R})^{\sigma}=\mathbb{C}$, and so the statement is
equivalent to:
\[
\mathcal{K}(\mathbb{C})_{\#}\simeq\mathcal{K}(M_{2}(\mathbb{R}))_{\#}%
^{h\sigma}\text{,}%
\]
an analogue of the well-known fact in Lemma \ref{hZ/2 for C} that
\[
\mathcal{K}(\mathbb{R})\simeq\mathcal{K}(\mathbb{C})^{\tau}%
\]
where $\tau$ denotes complex conjugation. The statement (\ref{K(M2(R)) fixed
points}) is a particular case of the following theorem.

\begin{Theorem}
Let $A$ be a real or complex simple algebra, and let $\sigma$ be a nontrivial
involution. Then
\[
\omega:\mathcal{K}(A^{\sigma})_{\#}\longrightarrow(\mathcal{K}(A)^{h\sigma
})_{\#}%
\]
is a homotopy equivalence.
\end{Theorem}

\noindent\textbf{Proof.}\quad Let us consider the complex case first, which
means that $A=M_{n}(\mathbb{C})$.\textsf{ }By the Skolem-Noether theorem,
$\sigma$ is the inner conjugation by a matrix whose square is a scalar matrix.
So, after extracting square roots and effecting a change of basis of
$\mathbb{C}^{n}$, we can assume that the matrix is of type $\left(
\begin{array}
[c]{cc}%
I_{\ell} & 0\\
0 & -I_{q}%
\end{array}
\right)  $ where $0<\ell=n-q<n$, since $\sigma$ is nontrivial. Since we are
stabilizing and suspending, we are reduced to proving the $2$-adic equivalence%
\[
\mathcal{K}(\mathbb{C})\times\mathcal{K}(\mathbb{C})\longrightarrow
\mathcal{K}(M_{2}(\mathbb{C}))^{h\mathbb{Z}/2}%
\]
shown in Lemma \ref{(E,F)}.

When $A$ is a ring of matrices over $\mathbb{R}$ or $\mathbb{H}$, we consider
\[
A^{\prime}=A\otimes_{\mathbb{R}}\mathbb{C\cong M}_{n}(\mathbb{C}).
\]
Let $\tau$ denote the complex conjugation on $A^{\prime}$, and write $\sigma$
also to denote the extension of $\sigma$ to $A^{\prime}$. From the complex
case, we already have the homotopy equivalence%
\[
\mathcal{K}(A^{\prime\sigma})_{\#}\simeq(\mathcal{K}(A^{\prime})^{h\sigma
})_{\#}\text{.}%
\]
So the result follows from the commutative diagram%
\[%
\begin{array}
[c]{ccccc}%
\mathcal{K}(A^{\sigma}) & \overset{\simeq}{\longrightarrow} & \mathcal{K}%
((A^{\sigma})^{\prime})^{h\tau} & = & \mathcal{K}((A^{\prime})^{\sigma
})^{h\tau}\\
\downarrow^{\omega} &  &  &  & \downarrow\\
\mathcal{K}(A)^{h\sigma} & \overset{\simeq}{\longrightarrow} & \mathcal{K}%
((A^{\prime})^{h\tau})^{h\sigma} & = & \mathcal{K}((A^{\prime})^{h\sigma
})^{h\tau}%
\end{array}
\]
in which the horizontal arrows are the descent homotopy equivalences of
\cite{K: descent}.\hfill$\Box\smallskip$

Observe that in the above proof, the right-hand vertical arrow requires
$2$-adic completion to be an equivalence. Therefore this property is also true
of the map $\omega$.

$\smallskip$\
%
%
%
%
%
\newpage

\begin{Void}
\textbf{Relation to the integral orthogonal group conjecture.}
\end{Void}%

\nopagebreak
Consider the diagram%
\[%
\begin{array}
[c]{ccc}%
B{}_{\varepsilon}O(\mathbb{Z}^{\prime})_{\#}^{+} & \overset{\alpha^{\prime}%
}{\longrightarrow} & B\mathrm{GL}(\mathbb{Z}^{\prime})_{\#}^{+h({}%
_{\varepsilon}\mathbb{Z}/2)}\\
\vspace{-0.13in} & \  & \\
\uparrow^{\lambda_{O}} &  & \uparrow^{\lambda_{\mathrm{GL}}}\\
B{}_{\varepsilon}O(\mathbb{Z})_{\#}^{+} & \overset{\alpha}{\longrightarrow} &
B\mathrm{GL}(\mathbb{Z})_{\#}^{+h({}_{\varepsilon}\mathbb{Z}/2)}\\
\vspace{-0.13in} & \  & \\
\downarrow^{\rho_{O}} &  & \downarrow^{\rho_{\mathrm{GL}}}\\
B{}_{\varepsilon}O(\mathbb{F}_{p})_{\#}^{+} & \overset{\alpha_{p}%
}{\longrightarrow} & B\mathrm{GL}(\mathbb{F}_{p})_{\#}^{+h({}_{\varepsilon
}\mathbb{Z}/2)}%
\end{array}
\]
and its effect on homotopy groups of dimension at least $2$. From Quillen's
localization theorem, we know that $\lambda_{\mathrm{GL}}$ induces
isomorphisms in this range. Likewise, isomorphisms for $\alpha^{\prime}$ and
$\alpha_{p}$ are proved in Theorem \ref{Thm C in 7} and Lemma \ref{hZ/2 for
F_3} respectively. Thus our conjecture that $\lambda_{O}$ induces isomorphisms
is equivalent to the `Lichtenbaum-Quillen' conjecture that $\alpha$ does also,
again for the same range of dimension.

\bigskip

\begin{center}
\textsc{Appendix A: on suspensions}\smallskip\label{suspension}
\end{center}%

\nopagebreak
For a \emph{discrete }ring $\Lambda$, we define the cone $C\Lambda$ of
$\Lambda$ to be the ring of infinite matrices (indexed by $\mathbb{N}$) over
$\Lambda$ for which there exists a natural number that bounds:

\begin{enumerate}
\item [(i)]the number of nonzero entries in each row and column; and

\item[(ii)] the number of distinct entries in the entire matrix.
\end{enumerate}

If $\tilde{\Lambda}$ denotes the ideal of $C\Lambda$ comprising matrices with
only finitely many nonzero entries, then the suspension $S\Lambda$ is defined
to be the quotient ring $C\Lambda/\tilde{\Lambda}$. Note that we have
canonical isomorphisms
\[
C\Lambda=\Lambda\otimes_{\mathbb{Z}}C\mathbb{Z},\qquad\qquad S\Lambda
=\Lambda\otimes_{\mathbb{Z}}S\mathbb{Z}\text{;}%
\]
from the second we readily obtain
\[
S(\Lambda\lbrack x])=(S\Lambda)[x]\text{.}%
\]

On the other hand, if $\Lambda$ is $\mathbb{R}$, $\mathbb{C}$ or\textsf{
}$\mathbb{H}$, with the usual topology, then we define $C\Lambda$ to be the
ring of bounded operators in a real, complex or quaternionic Hilbert space $H$
(with an $\ell^{2}$-basis) and $\tilde{\Lambda}$ to be the ideal of compact
operators in $H$. Then $S\Lambda$, defined to be $C\Lambda/\tilde{\Lambda}$,
is called the Calkin algebra. It is related to the space of Fredholm operators
$\mathcal{F}(H)$ by the surjective homotopy equivalence%
\[
\mathcal{F}(H)\longrightarrow\mathrm{GL}_{1}(S\Lambda)\text{.}%
\]
By virtue of Kuiper's theorem on the contractibility of $\mathrm{Aut}(H)$,
there are also chains of homotopy equivalences%
\[
\mathrm{GL}_{1}(S\Lambda)\longrightarrow\mathrm{GL}_{2}(S\Lambda
)\longrightarrow\cdots\longrightarrow\mathrm{GL}(S\Lambda)
\]
and%
\[
\mathcal{F}(H)\longrightarrow\mathcal{F}(H\oplus H)\longrightarrow
\cdots\longrightarrow\mathcal{F}(H\oplus\cdots\oplus H\oplus\cdots)
\]
as discussed in \cite{Karoubi: Espaces classifiants}.

With these definitions, there are obvious maps
\[
S\mathbb{Z}^{\prime}\rightarrow S\mathbb{R}\rightarrow S\mathbb{C}\rightarrow
S\mathbb{H}%
\]
inducing%
\[
\mathcal{K}(\mathbb{Z}^{\prime})\rightarrow\mathcal{K}(\mathbb{R}%
)\rightarrow\mathcal{K}(\mathbb{C})\rightarrow\mathcal{K}(\mathbb{H})
\]
and (with the trivial antiinvolutions on $\mathbb{Z}^{\prime}$, $\mathbb{R}$
and $\mathbb{C}$)%
\[
{}_{\varepsilon}\mathcal{L}(\mathbb{Z}^{\prime})\rightarrow{}_{\varepsilon
}\mathcal{L}(\mathbb{R})\rightarrow{}_{\varepsilon}\mathcal{L}(\mathbb{C}%
)\text{.}%
\]
See Section \ref{main results}.

\bigskip

\begin{center}
\textsc{Appendix B: on equalizers}\smallskip\label{equalizer}
\end{center}

The purpose of this appendix is to review some background material that we
need for analysis of the homotopy fiber on $\mathcal{K}$ and $\mathcal{L}$
spaces of the difference map known as $\psi-1$. Our requirement is that this
map $\psi$ be induced by a map of rings. So, as in Theorem \ref{M-V}, to
obtain a Mayer-Vietoris sequence, first we need a method for replacing an
arbitrary ring homomorphism by an epimorphism.

\smallskip

\textbf{Mapping cylinder. }Given any homomorphism $\theta:A\rightarrow B$ of
$\mathbb{Z}[1/2\ell]$-algebras, construct the pull-back
\[
M(\theta)=\{(a,\ell(x))\in A\times B[x]\mid\theta(a)=\ell(0)\},
\]
and so the diagram
\[%
\begin{array}
[c]{ccc}%
M(\theta) & \longrightarrow &  B[x]\\
\quad\downarrow^{e_{0}^{\prime}} & \searrow^{\tilde{\theta}} & \quad
\downarrow^{e_{0}}\\
A & \overset{\theta}{\longrightarrow} & B
\end{array}
\]
Here $\tilde{\theta}:(a,\,\ell(x))\mapsto\ell(1)$, $e_{0}:\ell(x)\mapsto
\ell(0)$, and $e_{0}^{\prime}$, are epimorphisms. (Note that $b\in B$ is the
image under $\tilde{\theta}$ of $(1,\,1+(b-1)x)\in M(\theta)$.) The triangles
don't commute. However, because $\ell$ is invertible in $B$, homotopy
invariance applies (for mod\ $\ell$ $K$-groups $\tilde{K}$ by \cite{W}, and
for mod\ $\ell$ $L$-groups $\tilde{L}$ by Theorem \ref{htpy invar} above);
thus $e_{0}$ and $e_{1}$ induce the same isomorphism. Also, $f:A\rightarrow
M(\theta)$ sending $a$ to $(a,\,\theta(a\hat{)})$, where $\theta(a\hat{)}$
denotes the constant polynomial at $\theta(a)$, has $\tilde{\theta}\circ
f=\theta$ and $e_{0}^{\prime}\circ f=\mathrm{id}$. Then the Mayer-Vietoris
sequence (for mod\ $\ell$ $K$-groups $\tilde{K}$ by \cite{W}, and for
mod\ $\ell$ $L$-groups $\tilde{L}$ by Theorem \ref{M-V} above) makes the
induced homomorphism $f_{\ast}$ an isomorphism.

\smallskip

\textbf{Equalizer. }Suppose that we are given homomorphisms $f_{1}%
,f_{2}:A\rightarrow B$ of $\mathbb{Z}[1/2\ell]$-algebras. Consider the ring
$\bar{B}=B\times B$, with diagonal homomorphism $\Delta:B\rightarrow\bar{B}$.
Use the mapping cylinder to replace $\Delta$ by $\tilde{\Delta}:M(\Delta
)\rightarrow\bar{B}$. Now apply to%
\[%
\begin{array}
[c]{ccc}%
F & \overset{}{\twoheadrightarrow} & A\\
\downarrow^{\gamma} & ^{\cdot}\lrcorner & \quad\downarrow^{(f_{1},f_{2})}\\
M(\Delta) & \overset{\tilde{\Delta}}{\twoheadrightarrow} & \bar{B}%
\end{array}
\]
to obtain its exact Mayer-Vietoris sequence as the lower sequence of {\small
\[%
\begin{array}
[c]{ccccccccc}%
\cdots\rightarrow & 0 & \rightarrow & \tilde{K}_{n}(M(\Delta)) & \rightarrow &
\tilde{K}_{n}B & \rightarrow & 0 & \rightarrow\cdots\\
& \downarrow &  & \downarrow &  & \downarrow^{\Delta_{\ast}} &  & \downarrow &
\\
\cdots\rightarrow & \tilde{K}_{n}F & \rightarrow & \tilde{K}_{n}A\oplus
\tilde{K}_{n}(M(\Delta)) & \rightarrow & \tilde{K}_{n}B\oplus\tilde{K}_{n}B &
\rightarrow & \tilde{K}_{n-1}F & \rightarrow\cdots
\end{array}
\]
} On passing to the cokernel of the vertical injection map, we retrieve the
desired exact sequence%
\[
\cdots\rightarrow\tilde{K}_{n}F\rightarrow\tilde{K}_{n}A\overset{f_{1\ast
}-f_{2\ast}}{\longrightarrow}\tilde{K}_{n}B\rightarrow\tilde{K}_{n-1}%
F\rightarrow\cdots
\]
Similarly for $\tilde{L}$, making $F$ the ring we seek with $\mathcal{\tilde
{K}}(F)$ and $\mathcal{\tilde{L}}(F)$ the desired homotopy fiber.

\end{document}